\documentclass[11pt, reqno]{amsart}
\usepackage{ amsmath, mathrsfs, amsthm, amscd, amsfonts, amssymb, graphicx, color, enumerate}
\usepackage[bookmarks=false, bookmarksnumbered, colorlinks, plainpages, backref]{hyperref}
\definecolor{TITLE}{rgb}{0.0,0.0,1.0}
\definecolor{AUTHOR1}{rgb}{0.00,0.59,0.00}
\definecolor{AUTHOR2}{rgb}{0.50,0.00,1.00}
\definecolor{SECTION}{rgb}{0.50,0.00,1.00}
\definecolor{FOOTTITLE}{rgb}{0.00,0.50,0.75}
\definecolor{THM}{rgb}{0.7,0.3,0.3}
\definecolor{SEC}{rgb}{0.6,0.1,.5}

\makeatletter
\@namedef{subjclassname@2010}{%
	\textup{2010} Mathematics Subject Classification}
\makeatother
\newtheorem{theorem}{{\color{THM} Theorem}}[section]

\newtheorem{lemma}[theorem]{{\color{THM}Lemma}}
\newtheorem{proposition}[theorem]{{\color{THM}Proposition}}
\newtheorem{corollary}[theorem]{{\color{THM}Corollary}}
\theoremstyle{definition}
\newtheorem{definition}[theorem]{{\color{THM}Definition\ }}
\newtheorem{example}[theorem]{{\color{THM}Example}}

\newtheorem{remark}[theorem]{{\color{THM}Remark}}
\numberwithin{equation}{section}

\begin{document}
	\title[Noncommutative weak almost periodicity and left uniform continuity...]{Noncommutative weak almost periodicity  \\ and left uniform continuity over locally compact quantum groups}
	\author[R. Faal, A. A. Khosravi, and  H.  R. Ebrahimi Vishki]{Ramin Faal$^1$,  Amin A. Khosravi$^2$, and Hamid Reza Ebrahimi Vishki$^3$}

\address{$^1$ Department of Pure Mathematics, Ferdowsi University
of Mashhad, P. O. Box 1159, Mashhad 91775, Iran.}
\email{\textcolor[rgb]{0.00,0.00,0.84}{faal.ramin@yahoo.com}}
\curraddr{$^1$ Water, Energy, and Environmental Engineering
Research Unit, University of Oulu,
P.O. Box 4300, FI-90014 Oulu, Finland.}
\email{\textcolor[rgb]{0.00,0.00,0.84}{ramin.faal@oulu.fi}}

	\address{$^2$ Department of Pure Mathematics, Ferdowsi University of Mashhad, P.O. Box 1159, Mashhad 91775, IRAN.}
	\email{\textcolor[rgb]{0.00,0.00,0.84}{amin.khosravi@mail.um.ac.ir}}
	\address{$^3$  Department of Pure Mathematics,  Centre of Excellence
		in Analysis on Algebraic Structures (CEAAS), Ferdowsi University of
		Mashhad, P.O. Box 1159, Mashhad 91775, IRAN.}
	\email{\textcolor[rgb]{0.00,0.00,0.84}{vishki@um.ac.ir}}

	

\begin{abstract}
We study weak almost periodicity and left uniform continuity for
Hopf-von Neumann algebras and locally compact quantum groups.
For a Hopf-von Neumann algebra $(M,\Gamma)$, we introduce the space
${\rm WAP}_{{\rm Bessel},\pi_u}(M_*)$ through norm convergence
conditions on rows and columns of slices associated with Bessel
families in the universal representation. We prove that this space
is the largest $C^*$-subalgebra of ${\rm WAP}(M_*)$ when separately continuous lift of the co-multiplication is multiplicative, and that the two
spaces coincide when $M$ is commutative.  Furthermore, the
algebra
$\{x\in C_b(\mathbb{G}):
\Gamma(x)\in M(C_0(\mathbb{G})\otimes C_b(\mathbb{G}))\}$
is characterized as the greatest left introverted $C^*$-subalgebra
of ${\rm LUC}(L^1(\mathbb{G}))$ satisfying a left Bessel condition
on $L^2(\mathbb{G})$.  We also show that ${\rm LUC}(L^1(\mathbb{G}))$ is a
$C^*$-algebra whenever the center of $C_0(\mathbb{G})$ contains an
element on which the bounded counit is non-zero.
\end{abstract}
\maketitle

\section{Introduction}

Let $(M,\Gamma)$ be a Hopf-von Neumann algebra with predual
$M_*$. The space ${\rm WAP}(M_*)$ of weakly almost periodic
functionals on $M_*$ consists of those elements of $M$ whose
associated translation operators are weakly compact. This
definition extends the classical notion of weak almost
periodicity on locally compact groups and provides a common
framework for studying weak almost periodicity in group
algebras, Fourier algebras, and quantum group convolution
algebras. A basic question concerns the compatibility of weak
almost periodicity with multiplication in $M$. If $M$ is
commutative or subhomogeneous, then ${\rm WAP}(M_*)$ is a $C^*$-subalgebra of $M$ \cite{Dawe, Run1}. In the general case,
this question leads naturally to the study of conditions
ensuring that products of weakly almost periodic elements
remain weakly almost periodic.

The co-multiplication provides a useful way to approach this
problem, since the translation operators associated with an
element $x\in M$ are obtained by slicing $\Gamma(x)$. However,
the factorization of slices of a product involves infinite
rows and columns, and their norm convergence requires
additional control. In this paper, we formulate such control
using Bessel families in Hilbert space. We introduce the
coefficient spaces $Q_c^H(I,M^*)$ and $Q_r^H(I,M^*)$ and
identify them isometrically with the canonical column and row
spaces over $M^*$. These identifications connect the
factorization of vector functionals with the operator space
structure of $M^*$.

For the universal representation $\pi_u$ of $M$, the resulting
Bessel conditions define a subspace
${\rm WAP}_{{\rm Bessel},\pi_u}(M_*)$ of ${\rm WAP}(M_*)$.
We prove that this subspace is a unital $C^*$-subalgebra of
$M$. We also characterize its elements by weak compactness
of the column and row inflations of their left and right
translation operators. Thus the Bessel conditions admit both
a description in terms of vector functionals and a
description in terms of amplified translation operators.

A further characterization involves the separately
continuous lift of the co-multiplication. For
$x\in{\rm WAP}(M_*)$, let
$\Gamma_u(x)\in M^{**}\,\overline\otimes\,M^{**}$ 
denote this lift, obtained from the universal representation (see \cite[Proposition 3.7]{Danon}).
We show that $x\in{\rm WAP}_{{\rm Bessel},\pi_u}(M_*)$ if and
only if $x^*x,xx^*\in{\rm WAP}(M_*)$ and
$\Gamma_u(x^*x)=\Gamma_u(x)^*\Gamma_u(x),$ and
$\Gamma_u(xx^*)=\Gamma_u(x)\Gamma_u(x)^*.$
In particular,
${\rm WAP}_{{\rm Bessel},\pi_u}(M_*)$ is the greatest
$C^*$-subalgebra of ${\rm WAP}(M_*)$ on which $\Gamma_u$ is
multiplicative.

These identities also clarify the relation with the
noncommutative separate continuity studied by Daws
\cite{Danon}. We identify the Bessel conditions with the
tensor and multiplier conditions arising in that work.
We then construct a Hopf-von Neumann algebra and an element
$x$ such that
$x,x^*x,xx^*\in{\rm WAP}(M_*),$ but
$x\notin{\rm WAP}_{{\rm Bessel},\pi_u}(M_*).$
Consequently, weak almost periodicity of $x$ and its two
squares does not by itself imply the Bessel conditions.
The distinction between the lifts of the squares and the
products of the lifts is therefore essential. We discuss
its implications for the identification used in
\cite[Theorems 3.8 and 5.1]{Danon}.

When $M$ is commutative, we prove that
${\rm WAP}_{{\rm Bessel},\pi_u}(M_*)
=
{\rm WAP}(M_*).$
This gives another proof of the $C^*$-algebra property
established in \cite[Theorem 4]{Dawe}. We also study weakly
almost periodic elements of a co-amenable locally compact
quantum group $\mathbb G$ and obtain a characterization
in terms of quasi-multiplier conditions on $\Gamma(x)$.
This relates weak compactness of translation operators to
tensor product conditions involving $C_0(\mathbb G)$.

The second part of the paper concerns left uniformly
continuous elements. For a locally compact group $G$,
the classical space ${\rm LUC}(G)$ is a $C^*$-algebra.
In the quantum group setting,
${\rm LUC}(L^1(\mathbb G))$ is known to be a $C^*$-algebra
when $\mathbb G$ is semi-regular \cite{HuN}. For co-amenable
$\mathbb G$, coproduct characterizations are available
under regularity \cite{Salm}, and the $C^*$-algebra property
also holds when $C_0(\mathbb G)$ has a bounded approximate
identity contained in its center \cite{Runde3}.

We prove that the latter approximate identity hypothesis
can be replaced by a condition involving a single central
element. More precisely, if $\mathbb G$ is co-amenable,
with bounded counit $\epsilon$, and there exists
$a\in Z(C_0(\mathbb G))$ such that $\epsilon(a)\neq0$, then
${\rm LUC}(L^1(\mathbb G))$ is a $C^*$-algebra.
We further consider the algebra
\[
{\rm LUC}(\mathbb G)
=
\left\{
x\in C_b(\mathbb G):
\Gamma(x)\in
M\!\left(C_0(\mathbb G)\otimes C_b(\mathbb G)\right)
\right\}
\]
studied in \cite{Salm}. For co-amenable $\mathbb G$, we
characterize this algebra as the greatest left introverted
$C^*$-subalgebra of ${\rm LUC}(L^1(\mathbb G))$ satisfying
the left Bessel condition introduced below. This condition
is formulated using normal vector functionals on
$L^2(\mathbb G)$.

The paper is organized as follows. Section 2 collects the
preliminaries. Section 3 introduces the Bessel coefficient
spaces and establishes the algebraic properties and
characterizations of
${\rm WAP}_{{\rm Bessel},\pi_u}(M_*)$, together with the
counterexample discussed above. Section 4 treats the
commutative case. Section 5 gives a coproduct
characterization of weakly almost periodic elements of
co-amenable locally compact quantum groups. Finally,
Section 6 establishes the central-element criterion for
the $C^*$-algebra property of
${\rm LUC}(L^1(\mathbb G))$ and the characterization of
${\rm LUC}(\mathbb G)$ through the left Bessel condition.

\section{Preliminaries}

\subsection{Hopf-von Neumann algebras and quantum groups}

For a $C^*$-algebra $A$, we denote its multiplier algebra and
the space of quasi-multipliers by $M(A)$ and $QM(A)$, respectively.
We use $\otimes$ for the injective operator space tensor
product, which agrees with the minimal tensor product for
$C^*$-algebras, and $\overline{\otimes}$ for the von Neumann algebra
tensor product. We follow the standard operator space notation
in \cite{BleOp}.

Let $C$ be an operator space in a unital $C^*$-algebra $B$, with
$1_B\in C$. We regard $(A\otimes C)^{**}$ as a subspace of
$(A\otimes B)^{**}$ and put
\[
\begin{split}
M(A\otimes C)
&=\bigl\{x\in(A\otimes C)^{**}:
       (a\otimes1)x,\ x(a\otimes1)\in A\otimes C
       \text{ for all }a\in A\bigr\},\\
QM(A\otimes C)
&=\bigl\{x\in(A\otimes C)^{**}:
       (a\otimes1)x(a'\otimes1)\in A\otimes C
       \text{ for all }a,a'\in A\bigr\}.
\end{split}
\]
The products in these definitions are taken in
$(A\otimes B)^{**}$. When the first tensor factor, rather
than the second, is the unital operator space, we use the
corresponding definitions with the tensor factors interchanged.

For a Hilbert space $H$, let $B(H)$ and $K(H)$ denote the spaces
of bounded and compact operators on $H$, respectively. Inner
products are linear in the first variable. For $\xi,\eta\in H$,
we define $\omega_{\xi,\eta}\in B(H)_*$ by
\[
\omega_{\xi,\eta}(u)=\langle u\xi,\eta\rangle
\qquad(u\in B(H)).
\]
For a $C^*$-algebra $A$, we use the usual dual module actions
\[
(a\cdot\mu)(b)=\mu(ba),
\qquad
(\mu\cdot a)(b)=\mu(ab)
\qquad(a,b\in A,\ \mu\in A^*).
\]

Let $\pi_u:A\to B(H_u)$ be the universal representation of $A$.
Each $\mu\in A^*$ admits a normal extension
$\widetilde\mu\in B(H_u)_*$ satisfying
$\mu(a)=\widetilde\mu(\pi_u(a))$ for $a\in A$.
Applying Cohen's factorization theorem to the essential
$K(H_u)$-module $B(H_u)_*$, we obtain
$\widetilde\mu=\widetilde\mu'\cdot\theta$ for some
$\widetilde\mu'\in B(H_u)_*$ and $\theta\in K(H_u)$. Thus
\[
\mu(a)=\widetilde\mu'\bigl(\theta\pi_u(a)\bigr)
\qquad(a\in A).
\]
In particular, multiplication by $\theta$ is performed on the
chosen normal extension before restriction to $A$.
When $A=M$ is a von Neumann algebra, the universal representation
is taken for its underlying $C^*$-algebra, so that
$\pi_u(M)''$ is canonically isomorphic to $M^{**}$; it is not
assumed to be a normal representation of $M$ with its given
predual.
For a set $I$, we write $F\Subset I$ to mean that $F$ is a finite
subset of $I$. Nets indexed by finite subsets are directed by
inclusion.

Let $M$ be a von Neumann algebra. A co-multiplication on $M$ is a
normal unital $*$-homomorphism
$\Gamma:M\to M\overline{\otimes} M$ satisfying
\[
(\Gamma\otimes{\rm id})\circ\Gamma
=({\rm id}\otimes\Gamma)\circ\Gamma.
\]
The pair $(M,\Gamma)$ is called a Hopf-von Neumann algebra.
Well-known examples are $L^\infty(G)$ and ${\rm VN}(G)$ for
a locally compact group $G$, with co-multiplications given by
\[
\begin{aligned}
\Gamma_a(f)(s,t)&=f(st)
&&\bigl(f\in L^\infty(G)\bigr),\\
\Gamma_s(\lambda(x))&=\lambda(x)\otimes\lambda(x)
&&\bigl(x\in G\bigr),
\end{aligned}
\]
where the first identity is understood almost everywhere and
$\lambda$ is the left regular representation of $G$ on $L^2(G)$.
For further information, see \cite{KV1,T}.

The co-multiplication induces an associative, completely
contractive product $\star$ on $M_*$, determined by
\begin{equation}
\langle x,\rho\star\tau\rangle
=\langle\Gamma(x),\rho\otimes\tau\rangle
\qquad(x\in M,\ \rho,\tau\in M_*).
\label{Prod}
\end{equation}
For $(L^\infty(G),\Gamma_a)$ this is convolution on $L^1(G)$,
whereas for $({\rm VN}(G),\Gamma_s)$ it is pointwise
multiplication on the Fourier algebra $A(G)$.

We also use $\star$ for the induced $M_*$-module operations on
$M$. For $x\in M$ and $\tau\in M_*$, these are
\[
\tau\star x=({\rm id}\otimes\tau)\Gamma(x),
\qquad
x\star\tau=(\tau\otimes{\rm id})\Gamma(x).
\]
Equivalently, for $a,b\in M_*$ and $\tau\in M$,
\[
\langle\tau\star a,b\rangle=\langle\tau,a\star b\rangle,
\qquad
\langle a\star\tau,b\rangle=\langle\tau,b\star a\rangle.
\]
For $m\in M^*$ and $\tau\in M$, define
$\tau\Diamond m,m\square\tau\in M$ by
\[
\langle\tau\Diamond m,a\rangle=\langle m,a\star\tau\rangle,
\qquad
\langle m\square\tau,a\rangle=\langle m,\tau\star a\rangle
\qquad(a\in M_*).
\]
The first and second Arens products on $M^*=(M_*)^{**}$ are
then given by
\begin{equation}
\langle m\square n,\tau\rangle=\langle m,n\square\tau\rangle,
\qquad
\langle m\Diamond n,\tau\rangle=\langle n,\tau\Diamond m\rangle
\quad(m,n\in M^*,\ \tau\in M).
\label{Arens}
\end{equation}
We denote the first topological center of $M^*$ by
\[
Z_t(M^*)=
\bigl\{m\in M^*:n\mapsto m\square n
\text{ is }\sigma(M^*,M)\text{-continuous}\bigr\}.
\]
For more details on Arens products, see \cite{Da}.

A norm-closed subspace $C\subseteq M$ satisfying
$C\star M_*\subseteq C$ is called left introverted if
$\nu\square x\in C$ for every $\nu\in C^*$ and $x\in C$, where
\[
\langle\nu\square x,\omega\rangle
=\langle\nu,x\star\omega\rangle
\qquad(\omega\in M_*).
\]
Right introversion is defined similarly using $\Diamond$ and
the condition $M_*\star C\subseteq C$. We call a subspace
introverted if it is both left and right introverted.

A Hopf-von Neumann algebra $(M,\Gamma)$ is called co-amenable
if $M_*$ has a bounded approximate identity. The algebra
$(L^\infty(G),\Gamma_a)$ is co-amenable for every locally compact
group $G$, while $({\rm VN}(G),\Gamma_s)$ is co-amenable if
and only if $G$ is amenable \cite{Lep}.

A locally compact quantum group is a quadruple
$(M,\Gamma,\phi,\psi)$, where $(M,\Gamma)$ is a
Hopf-von Neumann algebra and $\phi$ and $\psi$ are normal,
semifinite, faithful weights which are left and right invariant,
respectively. Thus
\[
\begin{aligned}
\phi\bigl((\omega\otimes{\rm id})\Gamma(x)\bigr)
&=\omega(1)\phi(x)
&&\bigl(\omega\in M_*^+,\ x\in M^+,\ \phi(x)<\infty\bigr),\\
\psi\bigl(({\rm id}\otimes\omega)\Gamma(x)\bigr)
&=\omega(1)\psi(x)
&&\bigl(\omega\in M_*^+,\ x\in M^+,\ \psi(x)<\infty\bigr).
\end{aligned}
\]
Writing
\[
M_\phi
=\operatorname{span}\{x\in M^+:\phi(x)<\infty\},
\qquad
M_\psi
=\operatorname{span}\{x\in M^+:\psi(x)<\infty\},
\]
we extend the weights linearly to these spaces. The invariance
identities then hold for arbitrary $\omega\in M_*$ and
$x\in M_\phi$ or $x\in M_\psi$, respectively.
There is a unitary antipode $R$, and we choose the right Haar
weight to be $\psi=\phi\circ R$. We call the quantum group
commutative if $M$ is commutative, and co-commutative if the
Banach algebra $M_*$ is commutative. For the general theory,
see \cite{KV1, KV2}.

Haar weights generalize Haar measures on locally compact groups.
For a locally compact group $G$, the quadruple
$G_a=(L^\infty(G),\Gamma_a,\phi_a,\psi_a)$ is a locally
compact quantum group, where the weights are integration against
left and right Haar measures. Similarly,
$G_s=({\rm VN}(G),\Gamma_s,\phi_s,\phi_s)$ is a
locally compact quantum group, where $\phi_s$ is the
Plancherel weight, which is both left and right invariant
\cite[Definition VII.3.2]{Tak}.

We denote a locally compact quantum group
$(M,\Gamma,\phi,\psi)$ by $\mathbb{G}$, and write
$L^\infty(\mathbb{G})=M$ and $L^1(\mathbb{G})=M_*$.
Let $L^2(\mathbb{G})$ be the GNS Hilbert space of $\phi$.
Using the standard identification with the GNS Hilbert space
of $\psi$, we regard $L^\infty(\mathbb{G})$ as a von Neumann subalgebra
of $B(L^2(\mathbb{G}))$. We denote by $V$ the right fundamental unitary
on $L^2(\mathbb{G})\otimes L^2(\mathbb{G})$. It satisfies the pentagonal relation
\begin{equation}
V_{12}V_{13}V_{23}=V_{23}V_{12}
\label{eq:pentagon}
\end{equation}
and implements the co-multiplication through
\[
\Gamma(x)=V(x\otimes1)V^*
\qquad\bigl(x\in L^\infty(\mathbb{G})\bigr).
\]
We use the usual leg-numbering notation; see \cite{KV1}.

We use the right regular realization of the dual quantum group
and write
\[
L^\infty(\widehat{\mathbb{G}})
=\overline{\bigl\{({\rm id}\otimes\omega)(V):
          \omega\in B(L^2(\mathbb{G}))_*\bigr\}}^{\,w^*}.
\]
This von Neumann algebra has a co-multiplication
\[
\widehat{{\Gamma}}(x)=\widehat V(x\otimes1)\widehat V^*,
\qquad \widehat V=\Sigma V^*\Sigma,
\]
where $\Sigma$ is the flip on $L^2(\mathbb{G})\otimes L^2(\mathbb{G})$.
With its left and right Haar weights, it is a locally compact
quantum group, denoted by $\widehat{\mathbb{G}}$. This realization is
unitarily equivalent to the usual left regular realization
of the dual; see \cite[Section 2]{HuN}. Under the canonical
identification, $\widehat{\widehat{\mathbb{G}}}=\mathbb{G}$.

As in the classical case, we put
\[
C_0(\mathbb{G})
=\overline{\bigl\{(\omega\otimes{\rm id})(V):
          \omega\in B(L^2(\mathbb{G}))_*\bigr\}}^{\,\|\cdot\|}.
\]
This is a $C^*$-algebra which is weak$^*$ dense in
$L^\infty(\mathbb{G})$, and the restriction of $\Gamma$ is a
nondegenerate $*$-homomorphism
\[
\Gamma:C_0(\mathbb{G})\longrightarrow
M\bigl(C_0(\mathbb{G})\otimes C_0(\mathbb{G})\bigr).
\]
We also denote its strict extension to $M(C_0(\mathbb{G}))$ by
$\Gamma$. We write $\mathcal{M}(\mathbb{G})=C_0(\mathbb{G})^*$ for the quantum measure
algebra. For the $C^*$-algebraic approach, see \cite{KV1}.

By Cohen's factorization theorem, each $\mu\in L^1(\mathbb{G})$ can
be written as $\mu=\mu'\cdot a$ for some
$\mu'\in L^1(\mathbb{G})$ and $a\in C_0(\mathbb{G})$. Likewise, each
$\mu\in \mathcal{M}(\mathbb{G})$ has such a factorization with
$\mu'\in \mathcal{M}(\mathbb{G})$ and $a\in C_0(\mathbb{G})$.
For factorization involving $C_0(\widehat{\mathbb{G}})$, we instead
factor a chosen normal extension to $B(L^2(\mathbb{G}))$.
Since $C_0(\widehat{\mathbb{G}})$ acts nondegenerately on $L^2(\mathbb{G})$,
for every $\mu\in L^1(\mathbb{G})$ there exist
$\widetilde\mu'\in B(L^2(\mathbb{G}))_*$ and
$\widehat a\in C_0(\widehat{\mathbb{G}})$ such that
\[
\mu(x)=\widetilde\mu'(\widehat a x)
\qquad\bigl(x\in L^\infty(\mathbb{G})\bigr).
\]
The analogous factorizations on the other side, or on both
sides, are understood in the same way.
With our realization of the dual, the unitary $V$ belongs to both
$
M\bigl(C_0(\widehat{\mathbb{G}})\otimes C_0(\mathbb{G})\bigr)$
and $
M\bigl(K(L^2(\mathbb{G}))\otimes C_0(\mathbb{G})\bigr).
$

\subsection{Weakly almost periodic functionals}

Let $(M,\Gamma)$ be a Hopf-von Neumann algebra. For $x\in M$,
define $L_x,R_x:M_*\to M$ by
\[
L_x(\tau)=x\star\tau=(\tau\otimes{\rm id})\Gamma(x),
\qquad
R_x(\tau)=\tau\star x=({\rm id}\otimes\tau)\Gamma(x).
\]
The space of weakly almost periodic functionals on $M_*$ is
\[
\begin{aligned}
{\rm WAP}(M_*)
&=\{x\in M:L_x\text{ is weakly compact}\}\\
&=\{x\in M:R_x\text{ is weakly compact}\}.
\end{aligned}
\]
It is a norm-closed operator system and an introverted
$M_*$-bimodule \cite[Proposition 5.7]{DL}. The two Arens products
induce the same product on ${\rm WAP}(M_*)^*$; we denote
this product by $\star$. With this product,
${\rm WAP}(M_*)^*$ is a Banach algebra whose multiplication
is separately weak$^*$ continuous; see \cite[Lemma 3.6 and Section 5]{Danon}.

For a locally compact group $G$, the identification
$M=L^\infty(G)$ gives
${\rm WAP}(M_*)={\rm WAP}(G)$, where the latter denotes
the classical space of weakly almost periodic functions,
viewed as a subspace of $L^\infty(G)$.

\subsection{The co-multiplication and the universal representation}

Let $\kappa_M:M\to M^{**}$ be the canonical embedding. The
universal representation $\pi_u:M\to B(H_u)$ extends to a normal
$*$-isomorphism from $M^{**}$ onto $\pi_u(M)''$. We use this
isomorphism to represent $M^{**}$ on $H_u$ and to identify
$\kappa_M(a)$ with $\pi_u(a)$ for $a\in M$.
As a $*$-homomorphism between $C^*$-algebras, $\Gamma$ has the
normal bidual extension
\[
\Gamma^{**}:M^{**}\longrightarrow
\bigl(M\overline{\otimes}M\bigr)^{**}.
\]
Its codomain is the bidual of $M\overline{\otimes}M$ and should
be distinguished from $M^{**}\overline{\otimes}M^{**}$.
In particular, the expression $(\pi_u\otimes\pi_u)\Gamma(x)$
cannot be interpreted as a normal tensor product of the original
representations, since $\pi_u$ need not be normal on $M$ with
its given predual.

For $x\in {\rm WAP}(M_*)$, there is a canonical realization of
the co-multiplication in $M^{**}\overline{\otimes}M^{**}$.
Indeed, weak compactness gives
$
L_x^{**}(M^*),R_x^{**}(M^*)\subseteq\kappa_M(M).
$
We identify these ranges with $M$ and define $\Gamma_u(x)$ by
\begin{equation}\label{universal coproduct pairing}
\langle\Gamma_u(x),\mu\otimes\nu\rangle
=\nu\bigl(L_x^{**}(\mu)\bigr)
=\mu\bigl(R_x^{**}(\nu)\bigr)
\qquad(\mu,\nu\in M^*).
\end{equation}
The two expressions agree because $L_x$ and $R_x$ are the
transpose orbit maps and are weakly compact. The bilinear form
in \eqref{universal coproduct pairing} is jointly completely
bounded and therefore determines an element of
$
M^{**}\overline{\otimes}M^{**}
=\bigl(M^*\widehat\otimes M^*\bigr)^*
$
where $\widehat\otimes$ denotes the operator space projective
tensor product. This is the separately continuous lift described
in \cite[Proposition 3.7]{Danon}. Its slices satisfy
\begin{equation}\label{universal coproduct slices}
\begin{aligned}
(\mu\otimes {\rm id})\Gamma_u(x)
&=\kappa_M\bigl(L_x^{**}(\mu)\bigr),\\
({\rm id}\otimes\nu)\Gamma_u(x)
&=\kappa_M\bigl(R_x^{**}(\nu)\bigr).
\end{aligned}
\end{equation}
 When $\mu,\nu\in M_*$, the pairing in
\eqref{universal coproduct pairing} is
$\langle\Gamma(x),\mu\otimes\nu\rangle$.
Thus $\Gamma_u(x)$ recovers the original co-multiplication when
tested against normal functionals, while also specifying its
slices by arbitrary elements of $M^*$.

Representing $M^{**}$ on $H_u$, we regard $\Gamma_u(x)$ as an
operator on $H_u\otimes H_u$. In this representation,
\[
(\omega_{\zeta,\eta}\otimes {\rm id})\Gamma_u(x)
=\pi_u\bigl(L_x^{**}(\omega_{\zeta,\eta})\bigr),
\]
where on the right
$\omega_{\zeta,\eta}(a)=\langle\pi_u(a)\zeta,\eta\rangle$.
The corresponding formula on the second tensor factor uses
$R_x^{**}$. We retain the notation $\Gamma_u$ to distinguish
this realization from both $\Gamma$ and $\Gamma^{**}$.
Multiplicativity of this lift on a specified subalgebra requires
a separate argument; it is not part of its definition on
${\rm WAP}(M_*)$.

\subsection{Left uniformly continuous elements}

Let $\mathbb{G}$ be a locally compact quantum group. The space of left
uniformly continuous elements is defined by
\[
{\rm LUC}(L^1(\mathbb{G}))
=\overline{\operatorname{span}
   \{x\star\omega:x\in L^\infty(\mathbb{G}),\ \omega\in L^1(\mathbb{G})\}}
   ^{\,\|\cdot\|}.
\]
If $\mathbb{G}$ is co-amenable, Cohen's factorization theorem gives
\[
{\rm LUC}(L^1(\mathbb{G}))
=L^\infty(\mathbb{G})\star L^1(\mathbb{G}),
\]
where the right-hand side denotes the set of individual
products $x\star\omega$.
As usual, we write $C_b(\mathbb{G})=M(C_0(\mathbb{G}))$. We have
\[
C_0(\mathbb{G})\subseteq{\rm LUC}(L^1(\mathbb{G}))\subseteq C_b(\mathbb{G}),
\]
and ${\rm LUC}(L^1(\mathbb{G}))$ is a left introverted operator
system \cite{HuN,Runde3}. For a classical locally compact group $G$,
this definition recovers ${\rm LUC}(G)$.

Salmi \cite{Salm} introduced the $C^*$-algebra
\[
{\rm LUC}(\mathbb{G})
=\bigl\{x\in C_b(\mathbb{G}):
\Gamma(x)\in
M\bigl(C_0(\mathbb{G})\otimes C_b(\mathbb{G})\bigr)\bigr\}.
\]
We retain the distinct notations ${\rm LUC}(\mathbb{G})$ and
${\rm LUC}(L^1(\mathbb{G}))$. If $\mathbb{G}$ is regular and co-amenable,
then
\[
{\rm LUC}(\mathbb{G})={\rm LUC}(L^1(\mathbb{G}))
\]
by \cite{Salm}. In Section 6, for co-amenable $\mathbb{G}$, we characterize
${\rm LUC}(\mathbb{G})$ as the greatest left introverted
$C^*$-subalgebra of ${\rm LUC}(L^1(\mathbb{G}))$ satisfying
the left Bessel condition.

\section{Bessel families and weakly almost periodic functionals}
\label{Bessel section}

Let $X$ be an operator space, and let $I,J$ be sets. We denote
by $M_{I,J}(X)$ the space of all $I\times J$ matrices whose
finite submatrices have uniformly bounded norm, with
\[
\|[x_{i,j}]\|
=\sup\bigl\{\|[x_{i,j}]_{i\in F,j\in E}\|:
                  F\Subset I,\ E\Subset J\bigr\}.
\]
The matrix norms are defined using
$M_n(M_{I,J}(X))=M_{I,J}(M_n(X))$ in the same way.
We denote by $K_{I,J}(X)$ the norm closure of the finitely
supported matrices in $M_{I,J}(X)$; see \cite{BleOp}.
For a $C^*$-algebra $M$, the column and row norms are
\[
\|[a_i]\|_{M_{I,1}(M)}^2
=\sup_{F\Subset I}\left\|\sum_{i\in F}a_i^*a_i\right\|,
\qquad
\|[a_i]\|_{M_{1,I}(M)}^2
=\sup_{F\Subset I}\left\|\sum_{i\in F}a_i a_i^*\right\|.
\]
Consequently, a column belongs to $K_{I,1}(M)$ precisely when
$\sum_i a_i^*a_i$ converges in norm, and a row belongs to
$K_{1,I}(M)$ precisely when $\sum_i a_i a_i^*$ converges in norm.
All sums over arbitrary index sets are nets of finite partial sums.

 A family $(\eta_i)_{i\in I}$ in Hilbert space $H$
is called a Bessel family if there exists $B\geq0$ such that
\[
\sum_{i\in I}|\langle\zeta,\eta_i\rangle|^2
\leq B\|\zeta\|^2
\qquad(\zeta\in H).
\]
We write
\[
\|(\eta_i)\|_{\rm Bes}
=\sup_{\|\zeta\|\leq1}
  \left(\sum_{i\in I}|\langle\zeta,\eta_i\rangle|^2\right)^{1/2}.
\]
Equivalently, the synthesis operator
\[
S_\eta:\ell^2(I)\longrightarrow H,
\qquad S_\eta((c_i))=\sum_{i\in I}c_i\eta_i,
\]
is bounded. In this case,
$\|S_\eta\|=\|(\eta_i)\|_{\rm Bes}$ and
$S_\eta^*\zeta=(\langle\zeta,\eta_i\rangle)_{i\in I}$.

Throughout this section, $\pi:M\to B(H)$ is a unital
$*$-representation of the underlying $C^*$-algebra of $M$
such that every $\mu\in M^*$ can be written as
\begin{equation}\label{vector representation hypothesis}
\mu(a)=\omega_{\zeta,\eta}(a)
=\langle\pi(a)\zeta,\eta\rangle
\qquad(a\in M)
\end{equation}
for some $\zeta,\eta\in H$. In particular, $\pi$ is faithful.
The universal representation $\pi_u$ on $H_u$ has this property.
The notation $\omega_{\zeta,\eta}$ will mean the restriction
through $\pi$ when its argument belongs to $M$.

\begin{definition}\label{Bessel coefficient spaces}
For a set $I$ and unital $*$-representation $\pi:M\to B(H)$ with property \eqref{vector representation hypothesis}, define
\[
\begin{aligned}
Q_c^H(I,M^*)
&=\bigl\{[\omega_{\zeta,\eta_i}]_{i\in I}:
       \zeta\in H,\ (\eta_i)_{i\in I}\text{ is Bessel}\bigr\},\\
Q_r^H(I,M^*)
&=\bigl\{[\omega_{\zeta_i,\eta}]_{i\in I}:
       (\zeta_i)_{i\in I}\text{ is Bessel},\ \eta\in H\bigr\}.
\end{aligned}
\]
For a column $[\mu_i]\in Q_c^H(I,M^*)$, put
\[
\|[\mu_i]\|_{Q_c^H}
=\inf\bigl\{\|\zeta\|\,\|(\eta_i)\|_{\rm Bes}:
       \mu_i=\omega_{\zeta,\eta_i}\text{ for every }i\in I\bigr\}.
\]
For a row $[\mu_i]\in Q_r^H(I,M^*)$, put
\[
\|[\mu_i]\|_{Q_r^H}
=\inf\bigl\{\|(\zeta_i)\|_{\rm Bes}\,\|\eta\|:
       \mu_i=\omega_{\zeta_i,\eta}\text{ for every }i\in I\bigr\}.
\]
The infima are taken over all such representations of the
entire family in the fixed representation $\pi$.
\end{definition}

\begin{remark}\label{cyclic_embed}
Let $\pi:M\to B(H)$ be a unital $*$-representation such that
every $\mu\in M^*$ can be written as
$\mu(a)=\langle\pi(a)\alpha,\beta\rangle$ for $a\in M$, and
for some $\alpha,\beta\in H$.
We show that every cyclic representation of $M$ embeds
as a subrepresentation of $\pi$. Let $\varphi\in M^*$ be positive
and choose $\alpha,\beta\in H$ such that
$\varphi(a)=\langle\pi(a)\alpha,\beta\rangle$ for all
 $a\in M$.
Let $(\rho_\varphi,H_\varphi,\Lambda_\varphi)$ be its GNS
construction, and put
$H_\alpha=\overline{\pi(M)\alpha}.$
Since $\pi$ is a $*$-representation, $H_\alpha$ reduces $\pi(M)$.
Define the densely defined operator
\[
T_0:\pi(M)\alpha\subseteq H_\alpha\longrightarrow H_\varphi,
\qquad
T_0\bigl(\pi(a)\alpha\bigr)=\Lambda_\varphi(a).
\]
This is well defined, since
$\|\Lambda_\varphi(a)\|^2
=\varphi(a^*a)
=\langle\pi(a)\alpha,\pi(a)\beta\rangle,$ 
and hence $\pi(a)\alpha=0$ implies $\Lambda_\varphi(a)=0$.
Let $P_\alpha$ be the orthogonal projection of $H$ onto
$H_\alpha$. For $a,b\in M$, we have
\[
\begin{aligned}
\langle T_0\pi(a)\alpha,\Lambda_\varphi(b)\rangle
&=\varphi(b^*a)\\
&=\langle\pi(a)\alpha,\pi(b)\beta\rangle\\
&=\langle\pi(a)\alpha,P_\alpha\pi(b)\beta\rangle.
\end{aligned}
\]
It follows that
$\Lambda_\varphi(M)\subseteq\operatorname{dom}(T_0^*),$ and 
$T_0^*\Lambda_\varphi(b)=P_\alpha\pi(b)\beta.$
Thus $T_0$ is a densely defined operator with a densely defined adjoint, and thus it is closable. Write $T=\overline{T_0}$. Its range
is dense in $H_\varphi$, since it contains
$\Lambda_\varphi(M)$.
For $a,b\in M$,
\[
T_0\pi(a)\pi(b)\alpha
=\Lambda_\varphi(ab)
=\rho_\varphi(a)T_0\pi(b)\alpha.
\]
Passing to the closure shows that $\pi(a)$ preserves
$\operatorname{dom}(T)$ and
$T\pi(a)h=\rho_\varphi(a)Th$ for every $h\in\operatorname{dom}(T)$.
In particular, for each unitary $v\in M$, $\pi(v)$ maps
$\operatorname{dom}(T)$ onto itself. Let
$T=V|T|$ 
be the polar decomposition of $T$. Applying uniqueness of polar
decomposition to
$T=\rho_\varphi(v)^*T\pi(v)|_{H_\alpha}$
gives
$V\pi(v)|_{H_\alpha}=\rho_\varphi(v)V.$ In fact, to justify this identity, put $U=\pi(v)|_{H_\alpha},$ and 
$W=\rho_\varphi(v).$
Both $U$ and $W$ are unitary. Since
$U\operatorname{dom}(T)=\operatorname{dom}(T)$ and
$TUh=WTh$ for every $h\in\operatorname{dom}(T)$, we have
$T=W^*TU$ as an equality of operators with the same domain.
Using $T=V|T|$, we obtain
\[
T=W^*V|T|U
  =(W^*VU)(U^*|T|U).
\]
We show that this is another polar decomposition of $T$.
The operator $U^*|T|U$ is positive self-adjoint, with domain
\[
U^*\operatorname{dom}(|T|)
=U^*\operatorname{dom}(T)
=\operatorname{dom}(T).
\]
Moreover, $W^*VU$ is a partial isometry, since
$(W^*VU)^*(W^*VU)=U^*V^*VU$
is an orthogonal projection. The kernel condition for the
canonical polar decomposition is also preserved as
\[
\ker(W^*VU)
=U^*(\ker V)
=U^*(\ker|T|)
=\ker(U^*|T|U).
\]
Thus uniqueness of the polar decomposition gives $V=W^*VU.$
Multiplying on the left by $W$ yields $VU=WV$, or equivalently,
\[
V\pi(v)|_{H_\alpha}=\rho_\varphi(v)V.
\]
Since the unitaries linearly span $M$, the same identity holds
for every $a\in M$. Moreover, the density of the range of $T$
gives $VV^*=1_{H_\varphi}$. Therefore
$V^*:H_\varphi\longrightarrow H_\alpha\subseteq H$ 
is an isometric intertwiner, with
\[
\pi(a)V^*=V^*\rho_\varphi(a)
\qquad(a\in M).
\]
Its range reduces $\pi(M)$. Since every cyclic representation
is unitarily equivalent to a GNS representation of a positive
functional, the required embedding property follows.
\end{remark}
\begin{remark}\label{Bessel matrix identification}
With the canonical dual operator
space structure on $M^*$, the identity maps on the coefficient
families give isometric identifications
\begin{equation}\label{Bessel column row isometries}
Q_c^H(I,M^*)=M_{I,1}(M^*),
\qquad
Q_r^H(I,M^*)=M_{1,I}(M^*).
\end{equation}
Let $[\mu_i]\in Q_c^H(I,M^*)$, with
$\mu_i=\omega_{\zeta,\eta_i}$. The associated column-valued map
$u:M\longrightarrow B(\mathbb C,\ell^2(I))$ 
is given by
\[
u(a)(1)=[\mu_i(a)]_{i\in I}=S_\eta^*\pi(a)\zeta.
\]
Since $\pi$ is completely contractive and the operator
$\lambda\mapsto\lambda\zeta$ has norm $\|\zeta\|$, we obtain
$\|u\|_{cb}\leq\|S_\eta\|\,\|\zeta\|.$ 
For each finite subset $F\subseteq I$, the definition of the
dual operator space matrix norms gives
\[
\|[\mu_i]_{i\in F}\|_{M_{F,1}(M^*)}
=
\bigl\|a\mapsto[\mu_i(a)]_{i\in F}\bigr\|_{cb}
\leq\|S_\eta\|\,\|\zeta\|.
\]
Taking the supremum over finite $F$, and then the infimum over
all such representations, yields
$\|[\mu_i]\|_{M_{I,1}(M^*)}
\leq\|[\mu_i]\|_{Q_c^H}.$

Conversely, let $[\mu_i]\in M_{I,1}(M^*)$ and put
$C=\|[\mu_i]\|_{M_{I,1}(M^*)}.$
For every $a\in M$ and finite $F\subseteq I$, we have
$\sum_{i\in F}|\mu_i(a)|^2\leq C^2\|a\|^2.$
Consequently, $[\mu_i(a)]_{i\in I}\in\ell^2(I)$, and the map
\[
u:M\longrightarrow B(\mathbb C,\ell^2(I)),
\qquad
u(a)(1)=[\mu_i(a)]_{i\in I},
\]
is well defined. Let $u_F$ denote its restriction to the output
coordinates in $F$. Taking finite coordinate projections at
each matrix level gives
\[
\|u\|_{cb}
=
\sup_{\substack{F\subseteq I\\ F\ {\rm finite}}}\|u_F\|_{cb}
=
\sup_{\substack{F\subseteq I\\ F\ {\rm finite}}}
\|[\mu_i]_{i\in F}\|_{M_{F,1}(M^*)}
=
C.
\]

The rectangular form of the completely bounded factorization
theorem gives a unital $*$-representation
$\rho:M\to B(K)$, a vector $\xi\in K$, and a bounded operator
$S:\ell^2(I)\to K$ such that
$u(a)(1)=S^*\rho(a)\xi,$ and 
$\|S\|\,\|\xi\|\leq\|u\|_{cb}$;
see \cite[Theorem~1.6]{PisierOS}. The rectangular form follows
by viewing $B(\mathbb C,\ell^2(I))$ as a corner of
$B(\mathbb C\oplus\ell^2(I))$.

Let
$K_0=\overline{\rho(M)\xi},$ 
and let $P$ be the orthogonal projection of $K$ onto $K_0$.
The subspace $K_0$ reduces $\rho(M)$ and contains $\xi$.
By the embedding property presented in Remark \ref{cyclic_embed}, there is an isometric
intertwiner $U:K_0\to H$ satisfying
$U\rho(a)|_{K_0}=\pi(a)U$ for 
 $a\in M$.
Writing $(\delta_i)_{i\in I}$ for the canonical basis of
$\ell^2(I)$, put
$\zeta=U\xi,$ and $\eta_i=UPS\delta_i.$
Then $(\eta_i)$ is Bessel, with synthesis operator $UPS$, and $\|(\eta_i)\|_{\rm Bes}=\|UPS\|\leq\|S\|.$
Furthermore,
\[
\begin{aligned}
\mu_i(a)
&=\langle S^*\rho(a)\xi,\delta_i\rangle\\
&=\langle\rho(a)\xi,S\delta_i\rangle\\
&=\langle\rho(a)\xi,PS\delta_i\rangle\\
&=\langle U\rho(a)\xi,UPS\delta_i\rangle\\
&=\langle\pi(a)\zeta,\eta_i\rangle.
\end{aligned}
\]
Hence $[\mu_i]\in Q_c^H(I,M^*)$ and
\[
\|[\mu_i]\|_{Q_c^H}
\leq\|\zeta\|\,\|UPS\|
\leq\|\xi\|\,\|S\|
\leq\|u\|_{cb}
=\|[\mu_i]\|_{M_{I,1}(M^*)}.
\]
This proves the column isometry.

For the row statement, consider
$\mu_i^*(a)=\overline{\mu_i(a^*)}$. 
If $v:M\to B(\ell^2(I),\mathbb C)$ is the row-valued map
associated with $[\mu_i]$, then the column-valued map associated
with $[\mu_i^*]$ is $a\mapsto v(a^*)^*$. Taking adjoints at
each matrix level shows that these maps have the same completely
bounded norm. Thus
$\|[\mu_i^*]\|_{M_{I,1}(M^*)}
=
\|[\mu_i]\|_{M_{1,I}(M^*)}.$
The identity
$\omega_{\zeta,\eta}^*=\omega_{\eta,\zeta}$ also holds.
Consequently, taking adjoints exchanges the row and column
Bessel representations without changing their factorization
costs. Applying the column result to $[\mu_i^*]$ proves the
row isometry.

The identifications in \eqref{Bessel column row isometries}
preserve the vector space operations on coefficient families.
In particular, the two definitions above give Banach spaces.
Their operator space structures may be taken to be those
inherited through these identifications.
Finally, let $(e_j)_{j\in J}$ be an orthonormal basis of $H$.
Then
$[\omega_{\zeta,e_j}]_{j\in J}\in Q_c^H(J,M^*),$ and $
[\omega_{e_j,\eta}]_{j\in J}\in Q_r^H(J,M^*),$
with
\[
\|[\omega_{\zeta,e_j}]_{j\in J}\|_{Q_c^H}
=\|\zeta\|,
\qquad
\|[\omega_{e_j,\eta}]_{j\in J}\|_{Q_r^H}
=\|\eta\|.
\]

\end{remark}

\begin{definition}\label{WAP Bessel definition}
Let $(M,\Gamma)$ be a Hopf-von Neumann algebra and let $\pi$
satisfy \eqref{vector representation hypothesis}. We define
${\rm WAP}_{{\rm Bessel},\pi}(M_*)$ to be the set of all
$x\in {\rm WAP}(M_*)$ such that, for every set $I$ and every
$[\mu_i]\in Q_c^H(I,M^*)$,
\begin{equation}\label{Bessel column condition}
\begin{aligned}
\bigl[L_x^{**}(\mu_i)\bigr]_{i\in I},
\bigl[R_x^{**}(\mu_i)\bigr]_{i\in I}&\in K_{I,1}(M)
\end{aligned}
\end{equation}
and, for every $[\nu_i]\in Q_r^H(I,M^*)$,
\begin{equation}\label{Bessel row condition}
\begin{aligned}
\bigl[L_x^{**}(\nu_i)\bigr]_{i\in I},
\bigl[R_x^{**}(\nu_i)\bigr]_{i\in I}&\in K_{1,I}(M)
\end{aligned}
\end{equation}
\end{definition}


In the following results, $(\pi_u,H_u)$ is the universal representation
of $M$. For $x\in {\rm WAP}(M_*)$, we use the separately continuous
lift $\Gamma_u(x)$, so that
\[
(\mu\otimes {\rm id})(\Gamma_u(x))=L_x^{**}(\mu),
\qquad
({\rm id}\otimes\mu)(\Gamma_u(x))=R_x^{**}(\mu)
\qquad(\mu\in M^*).
\]
The values of these maps are identified with elements of $M$.
The following theorem shows that the Bessel conditions on both
orbit maps are preserved under multiplication.

\begin{theorem}\label{cWAP algebra}
Let $(M,\Gamma)$ be a Hopf-von Neumann algebra. Then
${\rm WAP}_{{\rm Bessel},\pi_u}(M_*)$ is a unital  $C^*-$subalgebra of ${\rm WAP}(M_*)$.
\end{theorem}

\begin{proof}
Let $x,y\in {\rm WAP}_{{\rm Bessel},\pi_u}(M_*)$. We show that
$xy\in {\rm WAP}_{{\rm Bessel},\pi_u}(M_*)$.
Let $(e_\lambda)_{\lambda\in\Lambda}$ be an orthonormal basis
for $H_u$, and put
$
Z=\Gamma_u(x)\Gamma_u(y).$
We first establish the required row and column conditions for $Z$,
and then identify $Z$ with $\Gamma_u(xy)$.
Let $[\omega_{\zeta_\alpha,\eta}]_{\alpha\in I}
\in Q_r^{H_u}(I,M^*)$, where $(\zeta_\alpha)_{\alpha\in I}$
is Bessel and $\eta\in H_u$ is fixed. The usual matrix
multiplication identity gives
\begin{align}
\left[(\omega_{\zeta_\alpha,\eta}\otimes {\rm id})(Z)\right]_{\alpha\in I}
=
\underbrace{\left[
 (\omega_{e_\lambda,\eta}\otimes {\rm id})(\Gamma_u(x))
                       \right]_{\lambda\in\Lambda}}_{A}
\underbrace{\left[
 (\omega_{\zeta_\alpha,e_\lambda}\otimes {\rm id})(\Gamma_u(y))
                       \right]_{\lambda\in\Lambda,\alpha\in I}}_{B}.
\label{Bessel row factorization}
\end{align}
Since $x\in {\rm WAP}_{{\rm Bessel},\pi_u}(M_*)$, the first
matrix belongs to $K_{1,\Lambda}(M)$. Let
$S:\ell^2(I)\to H_u$ be the synthesis operator defined by
$S\delta_\alpha=\zeta_\alpha$. The second matrix represents
$\Gamma_u(y)(S\otimes1)$, after identifying $H_u$ with
$\ell^2(\Lambda)$. Hence
 $\|B\|\leq\|\Gamma_u(y)\|\,\|S\|$ and $B\in M_{\Lambda,I}(M).$
Moreover, for each $\lambda\in\Lambda$,
\[
\left[
 (\omega_{\zeta_\alpha,e_\lambda}\otimes {\rm id})(\Gamma_u(y))
\right]_{\alpha\in I}\in K_{1,I}(M),
\]
because $[\omega_{\zeta_\alpha,e_\lambda}]_{\alpha\in I}
\in Q_r^{H_u}(I,M^*)$.

For $F\Subset\Lambda$, let $A_F$ be obtained from $A$ by
retaining the entries indexed by $F$ and setting the others
equal to zero. Then $A_FB$ is a finite sum of rows belonging
to $K_{1,I}(M)$, and
\[
\|AB-A_FB\|\leq\|A-A_F\|\,\|B\|\longrightarrow0.
\]
It follows from \eqref{Bessel row factorization} that
$\left[(\omega_{\zeta_\alpha,\eta}\otimes {\rm id})(Z)
\right]_{\alpha\in I}\in K_{1,I}(M).$
Equivalently, the sum
\[
\sum_{\alpha\in I}
(\omega_{\zeta_\alpha,\eta}\otimes {\rm id})(Z)
(\omega_{\zeta_\alpha,\eta}\otimes {\rm id})(Z)^*
\]
is norm convergent.
Similarly, let $[\omega_{\zeta,\eta_\alpha}]_{\alpha\in I}
\in Q_c^{H_u}(I,M^*)$. In this case,
\begin{align}
\left[(\omega_{\zeta,\eta_\alpha}\otimes {\rm id})(Z)
\right]_{\alpha\in I}
\notag
=
\left[
 (\omega_{e_\lambda,\eta_\alpha}\otimes {\rm id})(\Gamma_u(x))
    \right]_{\alpha\in I,\lambda\in\Lambda}
\left[
 (\omega_{\zeta,e_\lambda}\otimes {\rm id})(\Gamma_u(y))
  \right]_{\lambda\in\Lambda}.
\label{Bessel column factorization}
\end{align}
The first matrix belongs to $M_{I,\Lambda}(M)$, by the
Bessel property of $(\eta_\alpha)_{\alpha\in I}$, and each
of its columns belongs to $K_{I,1}(M)$. The second matrix
belongs to $K_{\Lambda,1}(M)$. Truncating the second matrix
to finite subsets of $\Lambda$, as above, shows that
\[
\left[(\omega_{\zeta,\eta_\alpha}\otimes {\rm id})(Z)
\right]_{\alpha\in I}\in K_{I,1}(M).
\]
Repeating these factorizations on the second leg, using the
conditions on $R_x^{**}$ and $R_y^{**}$, gives the corresponding
row and column conditions for the slices
$({\rm id}\otimes\omega)(Z)$.

Taking $I$ to be a singleton, we see in particular that every
vector slice of $Z$, on either leg, belongs to $M$.
Since every functional in $M^*$ is a vector functional in the
universal representation, we obtain
$(\mu\otimes {\rm id})(Z),
({\rm id}\otimes\mu)(Z)\in M$ for each $(\mu\in M^*).$
Let $q:M^{**}\to M$ be the canonical normal $*$-homomorphism,
which is the identity on the canonical copy of $M$. Then
\[
(q\overline\otimes q)(Z)=\Gamma(x)\Gamma(y)=\Gamma(xy).
\]
Consequently, for $\omega\in M_*$,
$L_{xy}(\omega)=(\omega\otimes {\rm id})(Z),$
where the right-hand side is regarded as an element of $M$.

Let $(\omega_\alpha)$ be a net in $\operatorname{Ball}(M_*)$.
Since $\operatorname{Ball}(M^*)$ is weak$^*$ compact, after
passing to a subnet we may assume that
$\omega_\alpha\to\mu$ in $\sigma(M^*,M)$ for some $\mu\in M^*$.
For each $\nu\in M^*$, we have
\begin{align*}
\nu\bigl(L_{xy}(\omega_\alpha)\bigr)
&=\omega_\alpha\bigl(( {\rm id}\otimes\nu)(Z)\bigr)\\
&\longrightarrow
\mu\bigl(( {\rm id}\otimes\nu)(Z)\bigr)\\
&=\nu\bigl((\mu\otimes {\rm id})(Z)\bigr).
\end{align*}
Thus $L_{xy}(\omega_\alpha)$ converges weakly to an element
of $M$. Hence $L_{xy}$ is weakly compact and
$xy\in {\rm WAP}(M_*)$.
By uniqueness of the separately continuous lift,
\begin{equation}
\Gamma_u(xy)=\Gamma_u(x)\Gamma_u(y).
\label{Bessel product lift}
\end{equation}
The row and column conditions already proved for $Z$ now show
that $xy\in {\rm WAP}_{{\rm Bessel},\pi_u}(M_*)$.

${\rm WAP}_{{\rm Bessel},\pi_u}(M_*)$ is also self-adjoint.
Indeed, for $\mu\in M^*$ put
$\mu^\sharp(a)=\overline{\mu(a^*)}$. Then we have
$\omega_{\zeta,\eta}^\sharp=\omega_{\eta,\zeta},$ and
\[
L_{x^*}^{**}(\mu)=L_x^{**}(\mu^\sharp)^*, \quad
R_{x^*}^{**}(\mu)=R_x^{**}(\mu^\sharp)^*.
\]
The same identities on $M_*$ show that $x^*\in {\rm WAP}(M_*)$
whenever $x\in {\rm WAP}(M_*)$. Moreover, $\sharp$ interchanges
$Q_c^{H_u}(I,M^*)$ and $Q_r^{H_u}(I,M^*)$, while taking
adjoints interchanges $K_{I,1}(M)$ and $K_{1,I}(M)$.
Therefore, $x^*$ lies in ${\rm WAP}_{{\rm Bessel},\pi_u}(M_*)$
whenever $x\in {\rm WAP}_{{\rm Bessel},\pi_u}(M_*)$.
It follows that for every $x$ in ${\rm WAP}_{{\rm Bessel},\pi_u}(M_*)$ we get
$
x^*, x^*x, xx^*
\in {\rm WAP}_{{\rm Bessel},\pi_u}(M_*).
$
Also $\Gamma_u(x^*)=\Gamma_u(x)^*$.

We now show that ${\rm WAP}_{{\rm Bessel},\pi_u}(M_*)$
is norm closed. Let $(x_n)$ be a sequence in this space such
that $\|x_n-x\|\to0$ for some $x\in M$. Since
${\rm WAP}(M_*)$ is norm closed, we have
$x\in{\rm WAP}(M_*)$. Moreover, 
$\|L_x^{**}-L_y^{**}\|_{cb}\leq\|x-y\|,$ and 
 $\|R_x^{**}-R_y^{**}\|_{cb}\leq\|x-y\|$.
Fix an index set $I$ and
$[\mu_i]_{i\in I}\in Q_c^{H_u}(I,M^*)$. Using the
isometric identification
$Q_c^{H_u}(I,M^*)\cong M_{I,1}(M^*)$, we obtain
\[
\begin{aligned}
\big\|
[L_{x_n}^{**}(\mu_i)-L_x^{**}(\mu_i)]_{i\in I}
\big\|_{M_{I,1}(M)}
&\leq
\|L_{x_n}^{**}-L_x^{**}\|_{cb}
\|[\mu_i]_{i\in I}\|_{M_{I,1}(M^*)}
\\
&\leq
\|x_n-x\|\,
\|[\mu_i]_{i\in I}\|_{M_{I,1}(M^*)}
\longrightarrow0.
\end{aligned}
\]
For each $n$, the column
$[L_{x_n}^{**}(\mu_i)]_{i\in I}$ belongs to $K_{I,1}(M)$.
Since $K_{I,1}(M)$ is norm closed in $M_{I,1}(M)$, it follows
that
\[
[L_x^{**}(\mu_i)]_{i\in I}\in K_{I,1}(M).
\]
The same argument applies to $R_x^{**}$ and, using
$Q_r^{H_u}(I,M^*)\cong M_{1,I}(M^*)$, to the row conditions.
Hence $x\in{\rm WAP}_{{\rm Bessel},\pi_u}(M_*)$.
It is clear that this space contains $1_M$. 
Together with the linearity and closure under products and
adjoints established above, we proved that
${\rm WAP}_{{\rm Bessel},\pi_u}(M_*)$ is a unital
$C^*$-subalgebra of $M$ contained in ${\rm WAP}(M_*)$.

\end{proof}

\begin{corollary}
An element $x\in M$ belongs to
${\rm WAP}_{{\rm Bessel},\pi_u}(M_*)$ if and only if,
for every index set $I$, the four inflations
\[
\begin{aligned}
L_x^{I,1},\,R_x^{I,1}
&:K_{I,1}(M_*)\longrightarrow K_{I,1}(M),\\
L_x^{1,I},\,R_x^{1,I}
&:K_{1,I}(M_*)\longrightarrow K_{1,I}(M)
\end{aligned}
\]
are weakly compact.
\end{corollary}
\begin{proof}
Let $T=L_x$ or $T=R_x$. Since $T$ is completely bounded,
its row and column inflations are bounded. Under the canonical
identifications
\[
K_{I,1}(M_*)^{**}\cong M_{I,1}(M^*),
\qquad
K_{I,1}(M)^{**}\cong M_{I,1}(M^{**}),
\]
we have
$(T^{I,1})^{**}([\mu_i]_{i\in I})
=[T^{**}(\mu_i)]_{i\in I}.$
Hence, by Gantmacher's theorem, $T^{I,1}$ is weakly compact
if and only if
\[
[T^{**}(\mu_i)]_{i\in I}\in K_{I,1}(M)
\qquad
([\mu_i]_{i\in I}\in M_{I,1}(M^*)).
\]
Using the identification
$Q_c^{H_u}(I,M^*)\cong M_{I,1}(M^*)$,
this is precisely the column Bessel condition for $T^{**}$.
The analogous row identifications show that weak compactness
of $T^{1,I}$ is equivalent to the row Bessel condition.
Finally, taking $I$ to be a singleton shows that the stated
weak compactness conditions imply $x\in{\rm WAP}(M_*)$.
Applying the preceding equivalences to both $L_x$ and $R_x$,
for every index set $I$, proves the assertion.
\end{proof}
\begin{theorem}\label{op}
Let $(M,\Gamma)$ be a Hopf-von Neumann algebra and let
$x\in {\rm WAP}(M_*)$. The following conditions are equivalent:
\begin{enumerate}[(i)]
\item $x\in {\rm WAP}_{{\rm Bessel},\pi_u}(M_*)$.
\item $x^*x,xx^*\in {\rm WAP}(M_*)$ and
$\Gamma_u(x^*x)=\Gamma_u(x)^*\Gamma_u(x),
\Gamma_u(xx^*)=\Gamma_u(x)\Gamma_u(x)^*.$
\end{enumerate}
\end{theorem}

\begin{proof}
The implication $(i)\Rightarrow(ii)$ follows from
Theorem~\ref{cWAP algebra}, self-adjointness, and
\eqref{Bessel product lift}. Conversely, suppose that $(ii)$ holds. Let
$(e_\lambda)_{\lambda\in\Lambda}$ be an orthonormal basis for
$H_u$, and fix $\zeta\in H_u$. For $F\Subset\Lambda$, put
\[
\begin{aligned}
s_F&=\sum_{\lambda\in F}
 \bigl((\omega_{\zeta,e_\lambda}\otimes {\rm id})
                         (\Gamma_u(x))\bigr)^*
 (\omega_{\zeta,e_\lambda}\otimes {\rm id})(\Gamma_u(x)),\\
s&=(\omega_{\zeta,\zeta}\otimes {\rm id})(\Gamma_u(x^*x)).
\end{aligned}
\]
Since $x,x^*x\in {\rm WAP}(M_*)$, both $s_F$ and $s$
belong to $M$. By the identities in $(ii)$, we get
$0\leq s_F\leq s$.
Let $\mu\in M^*$ be positive. In the universal representation
we may write $\mu=\omega_{\xi,\xi}$ for some $\xi\in H_u$.
Then
\begin{align*}
\mu(s_F)
&=\sum_{\lambda\in F}
  \left\|(\omega_{\zeta,e_\lambda}\otimes {\rm id})
                              (\Gamma_u(x))\xi\right\|^2\\
&\longrightarrow
  \left\|\Gamma_u(x)(\zeta\otimes\xi)\right\|^2\\
&=\left\langle
  \Gamma_u(x)^*\Gamma_u(x)(\zeta\otimes\xi),
                                    \zeta\otimes\xi\right\rangle\\
&=\left\langle
  \Gamma_u(x^*x)(\zeta\otimes\xi),\zeta\otimes\xi
  \right\rangle
=\mu(s).
\end{align*}
Thus $(s_F)$ is an increasing net in $M^+$ which converges
to $s$ against every positive functional in $M^*$.
The functions $\mu\mapsto\mu(s-s_F)$ on the weak$^*$ compact
state space $S(M)$ are continuous and decrease pointwise to
zero. Dini's theorem therefore gives
\[
\|s-s_F\|=\sup_{\mu\in S(M)}\mu(s-s_F)\longrightarrow0.
\]
Consequently,
$\left[(\omega_{\zeta,e_\lambda}\otimes {\rm id})
                   (\Gamma_u(x))\right]_{\lambda\in\Lambda}
\in K_{\Lambda,1}(M).$
Applying the same argument to $\Gamma_u(x)^*$, using
$\Gamma_u(xx^*)=\Gamma_u(x)\Gamma_u(x)^*$, gives
$\left[(\omega_{e_\lambda,\eta}\otimes {\rm id})
                   (\Gamma_u(x))\right]_{\lambda\in\Lambda}
\in K_{1,\Lambda}(M)
$ for $\eta\in H_u.$
The same calculations on the second leg give the corresponding
basic row and column conditions for $R_x^{**}$.

Now let $[\omega_{\zeta,\eta_\alpha}]_{\alpha\in I}
\in Q_c^{H_u}(I,M^*)$, and let
$S:\ell^2(I)\to H_u$ be given by $S\delta_\alpha=\eta_\alpha$.
For $T=L_x^{**}$ or $R_x^{**}$, we have
\[
\left[T(\omega_{\zeta,\eta_\alpha})\right]_{\alpha\in I}
=S^*\left[T(\omega_{\zeta,e_\lambda})\right]_{\lambda\in\Lambda}
\in K_{I,1}(M).
\]
Similarly, if $[\omega_{\zeta_\alpha,\eta}]_{\alpha\in I}
\in Q_r^{H_u}(I,M^*)$ and $S\delta_\alpha=\zeta_\alpha$, then
\[
\left[T(\omega_{\zeta_\alpha,\eta})\right]_{\alpha\in I}
=\left[T(\omega_{e_\lambda,\eta})\right]_{\lambda\in\Lambda}S
\in K_{1,I}(M).
\]
Hence $x\in {\rm WAP}_{{\rm Bessel},\pi_u}(M_*)$.
\end{proof}

\begin{example}\label{Bessel converse counterexample}
We construct a Hopf-von Neumann algebra $(M,\Gamma)$ and
a self-adjoint element $x\in M$ such that
$x,x^*x,xx^*\in{\rm WAP}(M_*),$ and
$x\notin{\rm WAP}_{{\rm Bessel},\pi_u}(M_*).$
Let $K=\ell^2(\mathbb N)$, with its standard orthonormal
basis $(e_i)_{i\geq1}$, and put
\[
M=M_0\oplus M_1\oplus M_2\oplus M_\infty
=\mathbb C\oplus B(K)\oplus B(K\otimes K)\oplus\mathbb C.
\]
We regard $M\overline\otimes M$ as the direct sum of the
spaces $M_r\overline\otimes M_s$, where
$r,s\in\{0,1,2,\infty\}$. For
$a=(a_0,a_1,a_2,a_\infty)\in M$, define
\[
\Gamma(a)_{r,s}=
\begin{cases}
a_{r+s},
 &r,s\in\{0,1,2\},\ r+s\leq2,\\[2mm]
a_\infty(1_{M_r}\otimes1_{M_s}),
 &\text{otherwise}.
\end{cases}
\]
In the first case we use the canonical identifications
\[
M_{r+s}=M_r\overline\otimes M_s,
\qquad K^{\otimes0}=\mathbb C.
\]
Thus the indices record tensor powers up to order $2$;
whenever their sum exceeds $2$, or an index is $\infty$,
the scalar component $a_\infty$ is used.
Each component map is a normal unital $*$-homomorphism,
so the same is true of $\Gamma$. Moreover,
\[
\Gamma(a)_{0,s}=a_s
\qquad(s\in\{0,1,2,\infty\}),
\]
which proves injectivity. Under the canonical identifications
of the triple tensor products, both iterated coproducts have
$(r,s,t)$-component
\[
\begin{cases}
a_{r+s+t},
 &r,s,t\in\{0,1,2\},\ r+s+t\leq2,\\[2mm]
a_\infty
 (1_{M_r}\otimes1_{M_s}\otimes1_{M_t}),
 &\text{otherwise}.
\end{cases}
\]
Hence $\Gamma$ is coassociative, and $(M,\Gamma)$ is a
Hopf-von Neumann algebra.

Let $\Sigma\in B(K\otimes K)$ be the flip operator,
$\Sigma(\xi\otimes\eta)=\eta\otimes\xi,$
and set
\[
x=(0,0,\Sigma,0),
\qquad
y=(0,0,1_{K\otimes K},0).
\]
Since $\Sigma^*=\Sigma$ and $\Sigma^2=1_{K\otimes K}$,
we have
$x=x^*,$ and $x^*x=xx^*=x^2=y.$
The only nonzero components of $\Gamma(x)$ are
$\Gamma(x)_{0,2}=\Gamma(x)_{1,1}=\Gamma(x)_{2,0}=\Sigma.$
The corresponding
components of $\Gamma(y)$ are the identity operators.

We next show that $x,y\in{\rm WAP}(M_*)$. The predual of $M$ can be identified with
\[
M_*=
\mathbb C\oplus_1 S_1(K)\oplus_1 S_1(K\otimes K)
\oplus_1\mathbb C,
\]
where $S_1(K)$ are the trace class operators  on $K$. For $f=(\alpha,\rho,\sigma,\beta)\in M_*$, the preceding
description of the coproduct gives
\[
\begin{aligned}
L_x(f)
&=\bigl(\operatorname{Tr}(\sigma\Sigma),
         \rho,\alpha\Sigma,0\bigr),\\
L_y(f)
&=\bigl(\operatorname{Tr}(\sigma),
         \operatorname{Tr}(\rho)1_K,
         \alpha1_{K\otimes K},0\bigr).
\end{aligned}
\]
Here we have used
$(\rho\otimes{\rm id})(\Sigma)=\rho,$
where $\rho$ is viewed as a normal functional on the left
and as an operator on $K$ on the right. The same formulas
hold for $R_x$ and $R_y$.

The inclusion $S_1(K)\to B(K)$ factors as
\[
S_1(K)\longrightarrow S_2(K)\longrightarrow B(K),
\]
with $\|\rho\|\leq\|\rho\|_2\leq\|\rho\|_1$ where  $S_2(K)$ is the Hilbert-Schmidt operators on $K$. Since $S_2(K)$ is a Hilbert space, this inclusion is weakly
compact. The remaining terms in $L_x$ have finite rank,
while $L_y$ itself has finite rank. Hence $L_x$ and $L_y$
are weakly compact, and consequently
$x,x^*x,xx^*\in{\rm WAP}(M_*).$

We now exhibit a Bessel column for which the required
norm convergence fails. Define
\[
\mu_i(a)=\langle a_1e_1,e_i\rangle
\qquad(a\in M,\ i\geq1).
\]
The representation $a\mapsto a_1$ on $K$ is cyclic, with
cyclic vector $e_1$. It therefore occurs as a subrepresentation
of $\pi_u$. Explicitly, there is an isometry $U:K\to H_u$
such that
\[
\pi_u(a)U=Ua_1
\qquad(a\in M).
\]
Thus, with $\zeta=Ue_1$ and $\eta_i=Ue_i$,
\[
\mu_i(a)=\langle\pi_u(a)\zeta,\eta_i\rangle.
\]
The family $(\eta_i)_{i\geq1}$ is orthonormal, and hence Bessel.
Therefore
$[\mu_i]_{i\geq1}\in Q_c^{H_u}(\mathbb N,M^*).$

Let $E_{ij}\in B(K)$ be the matrix units,
$E_{ij}\xi=\langle\xi,e_j\rangle e_i.$
Each $\mu_i$ is normal and corresponds, under the above
predual identification, to $(0,E_{1i},0,0)$. Hence
\[
d_i:=L_x^{**}(\mu_i)
=L_x(\mu_i)
=(0,E_{1i},0,0).
\]
Since $E_{1i}^*E_{1i}=E_{ii}$, we obtain
\[
\sum_{i=1}^n d_i^*d_i
=
\left(0,\sum_{i=1}^nE_{ii},0,0\right).
\]
For every $m>n$,
\[
\left\|\sum_{i=n+1}^m d_i^*d_i\right\|
=
\left\|\sum_{i=n+1}^mE_{ii}\right\|
=1.
\]
Thus the partial sums are not norm Cauchy. Consequently,
$[d_i]_{i\geq1}\notin K_{\mathbb N,1}(M),$
so the column Bessel condition fails and
$x\notin{\rm WAP}_{{\rm Bessel},\pi_u}(M_*).$

Finally, $x^3=x$ and $x^2=y$, so
$C^*(1_M,x)=\operatorname{span}\{1_M,x,y\}.$
Since $1_M,x,y\in{\rm WAP}(M_*)$ and ${\rm WAP}(M_*)$
is a linear space, this finite-dimensional $C^*$-algebra
is contained in ${\rm WAP}(M_*)$. It is not contained in
${\rm WAP}_{{\rm Bessel},\pi_u}(M_*)$, since it contains $x$.
Thus the Bessel algebra need not contain every $C^*$-subalgebra
of $M$ contained in ${\rm WAP}(M_*)$.
\end{example}

We next compare our definition with the characterization of
separate continuity in \cite{Danon}. 
Recall that
\[
{\rm SC}(M\times M)
=
\left\{
X\in M^{**}\overline\otimes M^{**}:
(\mu\otimes{\rm id})(X),\,
({\rm id}\otimes\mu)(X)\in M
\quad(\mu\in M^*)
\right\}.
\]
Following \cite[Definition~3.3]{Danon}, write
\[
M\overset{\rm sc}{\otimes}M
=
\left\{
X\in{\rm SC}(M\times M):
X^*X,\ XX^*\in{\rm SC}(M\times M)
\right\}.
\]
For $x\in{\rm WAP}(M_*)$, let
$\Gamma_u(x)\in{\rm SC}(M\times M)$ 
be the separately continuous lift given by
\cite[Proposition~3.7]{Danon}. Thus
\[
\begin{aligned}
(\mu\otimes{\rm id})\Gamma_u(x)
&=L_x^{**}(\mu),\\
({\rm id}\otimes\mu)\Gamma_u(x)
&=R_x^{**}(\mu)
\qquad(\mu\in M^*).
\end{aligned}
\]

\begin{proposition}\label{Bessel Daws characterization}
Let $(M,\Gamma)$ be a Hopf-von Neumann algebra and let
$x\in{\rm WAP}(M_*)$. The following conditions are equivalent:
\begin{enumerate}
\item[(i)]
$x\in{\rm WAP}_{{\rm Bessel},\pi_u}(M_*)$;

\item[(ii)]
$\Gamma_u(x)\in
M\!\left(\mathcal K(H_u)\otimes M\right)
\cap
M\!\left(M\otimes\mathcal K(H_u)\right);$

\item[(iii)]
$\Gamma_u(x)\in M\overset{\rm sc}{\otimes}M.$
\end{enumerate}
Consequently,
\[
{\rm WAP}_{{\rm Bessel},\pi_u}(M_*)
=
\left\{
x\in{\rm WAP}(M_*):
\Gamma_u(x)\in M\overset{\rm sc}{\otimes}M
\right\}.
\]
\end{proposition}

\begin{proof}
The equivalence of {\rm(ii)} and {\rm(iii)} follows from
\cite[Theorem~3.1]{Danon}, applied to the universal
representation of $M$. We prove that {\rm(i)}
and {\rm(ii)} are equivalent.
Suppose that {\rm(i)} holds, and fix an orthonormal basis
$(e_i)_{i\in I}$ of $H_u$. For $\zeta\in H_u$, put
\[
a_i=(\omega_{\zeta,e_i}\otimes{\rm id})\Gamma_u(x),
\qquad
b_i=(\omega_{e_i,\zeta}\otimes{\rm id})\Gamma_u(x).
\]
Since
$[\omega_{\zeta,e_i}]_{i\in I}
\in Q_c^{H_u}(I,M^*),$ and 
$[\omega_{e_i,\zeta}]_{i\in I}
\in Q_r^{H_u}(I,M^*),$ 
the Bessel conditions give
\begin{equation}\label{K-incl}
    [a_i]_{i\in I}\in K_{I,1}(M),
\qquad
[b_i]_{i\in I}\in K_{1,I}(M).
\end{equation}
For $\zeta,\eta\in H_u$, let $\theta_{\zeta,\eta}(\xi)
=\langle\xi,\eta\rangle\zeta.$
The matrix expansions
\[
\begin{aligned}
\Gamma_u(x)(\theta_{\zeta,\eta}\otimes1)
&=\sum_{i\in I}\theta_{e_i,\eta}\otimes a_i,\\
(\theta_{\eta,\zeta}\otimes1)\Gamma_u(x)
&=\sum_{i\in I}\theta_{\eta,e_i}\otimes b_i
\end{aligned}
\]
converge in norm. Indeed, for every finite $E\subseteq I$,
\[
\begin{aligned}
\left\|
\sum_{i\in E}\theta_{e_i,\eta}\otimes a_i
\right\|
&=
\|\eta\|
\left\|\sum_{i\in E}a_i^*a_i\right\|^{1/2},\\
\left\|
\sum_{i\in E}\theta_{\eta,e_i}\otimes b_i
\right\|
&=
\|\eta\|
\left\|\sum_{i\in E}b_i b_i^*\right\|^{1/2}.
\end{aligned}
\]
The corresponding tail norms tend to zero by equation \eqref{K-incl}. Hence
\[
\Gamma_u(x)(\theta\otimes1),\,
(\theta\otimes1)\Gamma_u(x)
\in\mathcal K(H_u)\otimes M
\qquad(\theta\in\mathcal K(H_u)).
\]
Applying the same argument to $R_x^{**}$ gives
$\Gamma_u(x)\in M\!\left(M\otimes\mathcal K(H_u)\right).$
Thus {\rm(ii)} holds.

Conversely, suppose that {\rm(ii)} holds. Fix
$\zeta\in H_u$, choose a unit vector $\eta\in H_u$, and let
$P_F$ be the orthogonal projection onto
$\operatorname{span}\{e_i:i\in F\}$, where $F$ runs through
the finite subsets of $I$. Since
\[
\Gamma_u(x)(\theta_{\zeta,\eta}\otimes1)
\in\mathcal K(H_u)\otimes M,
\]
we have
\[
\begin{aligned}
\left\|
\left[
L_x^{**}(\omega_{\zeta,e_i})
\right]_{i\in I\setminus F}
\right\|
=
\left\|
((1-P_F)\otimes1)
\Gamma_u(x)(\theta_{\zeta,\eta}\otimes1)
\right\|\to 0.
\end{aligned}
\]
Therefore
$\left[
L_x^{**}(\omega_{\zeta,e_i})
\right]_{i\in I}
\in K_{I,1}(M).$ 
Applying the same argument to $\Gamma_u(x)^*$ and taking
adjoints gives
$\left[
L_x^{**}(\omega_{e_i,\zeta})
\right]_{i\in I}
\in K_{1,I}(M).$
The second multiplier condition gives the corresponding
column and row assertions for $R_x^{**}$.

It remains to pass from orthonormal families to arbitrary
Bessel families. Let $(\eta_j)_{j\in J}$ be Bessel, with
synthesis operator
\[
S_\eta:\ell^2(J)\longrightarrow H_u,
\qquad S_\eta\delta_j=\eta_j.
\]
Using $(e_i)_{i\in I}$ to identify $H_u$ with $\ell^2(I)$,
we obtain
\[
\left[
L_x^{**}(\omega_{\zeta,\eta_j})
\right]_{j\in J}
=
S_\eta^*
\left[
L_x^{**}(\omega_{\zeta,e_i})
\right]_{i\in I}
\in K_{J,1}(M).
\]
Similarly, if $(\zeta_j)_{j\in J}$ is Bessel, with synthesis
operator $S_\zeta$, then
\[
\left[
L_x^{**}(\omega_{\zeta_j,\eta})
\right]_{j\in J}
=
\left[
L_x^{**}(\omega_{e_i,\eta})
\right]_{i\in I}S_\zeta
\in K_{1,J}(M).
\]
The same identities hold with
$R_x^{**}$ in place of $L_x^{**}$. Thus all four Bessel
conditions hold, and {\rm(i)} was proved.
\end{proof}

\begin{remark}\label{Bessel Daws square distinction}
The preceding proposition identifies our Bessel conditions
with tensor and multiplier conditions in \cite{Danon}. These must be
distinguished from
\[
x,\ x^*x,\ xx^*\in{\rm WAP}(M_*),
\tag{$*$}
\]
which is the condition used to introduce
${\rm WAP}(M,\Gamma)$ in \cite[Theorem~5.1]{Danon}.
Condition $(*)$ ensures that
\[
\Gamma_u(x),
\Gamma_u(x^*x),
\Gamma_u(xx^*)
\in{\rm SC}(M\times M).
\]
In contrast, the tensor condition in the proposition requires
\[
\Gamma_u(x)^*\Gamma_u(x),\quad
\Gamma_u(x)\Gamma_u(x)^*
\in{\rm SC}(M\times M).
\]
The separate continuity of the lifted squares does not by
itself establish separate continuity of the products of
the lifts. Example~\ref{Bessel converse counterexample}
shows that $(*)$ does not imply the Bessel conditions.
More precisely, let
$q:M^{**}\longrightarrow M$
be the canonical normal $*$-homomorphism, and put
$\Theta=q\overline\otimes q$.
Suppose that $(*)$ holds, and write
\[
X=\Gamma_u(x),
\qquad
Y=\Gamma_u(x^*x).
\]
Then $X,Y\in{\rm SC}(M\times M)$ and
\[
\begin{aligned}
\Theta(Y)
&=\Gamma(x^*x)\\
&=\Gamma(x)^*\Gamma(x)\\
&=\Theta(X)^*\Theta(X)\\
&=\Theta(X^*X).
\end{aligned}
\]
By \cite[Proposition~3.7]{Danon}, the restriction of
$\Theta$ to ${\rm SC}(M\times M)$ is injective.
However, at this stage $X^*X$ has not been shown to belong
to ${\rm SC}(M\times M)$. Hence the equality
$\Theta(Y)=\Theta(X^*X)$ does not justify the conclusion
$Y=X^*X$.

This is the unjustified inference in the converse argument
of \cite[Theorem~3.8]{Danon}. The missing point is that
injectivity on ${\rm SC}(M\times M)$ can be applied to
$X^*X$ only after its membership in that space has been
established.
\end{remark}


\section{Weakly almost periodic functions on commutative Hopf von Neumann algebras}

In this section, we show that, for every commutative Hopf-von Neumann
algebra $(M,\Gamma)$,
${\rm WAP}_{{\rm Bessel},\pi_u}(M_*)={\rm WAP}(M_*).$
Together with Theorem~\ref{cWAP algebra}, this gives another proof
that ${\rm WAP}(M_*)$ is a $C^*$-algebra, a result established
in \cite{Dawe}. We first record the following consequence of
the dominated convergence theorem.

\begin{lemma}\label{DCT}
Let $X$ be a compact Hausdorff space, let $\mu$ be a finite
positive regular Borel measure on $X$, and let $\mathcal F$
be a weakly compact subset of $C(X)$. Then the canonical map
\[
\mathcal F\longrightarrow L^2(X,\mu),\qquad g\longmapsto g,
\]
is continuous from the weak topology of $\mathcal F$ to the
norm topology of $L^2(X,\mu)$.
\end{lemma}

\begin{proof}
Let $J:C(X)\to L^2(X,\mu)$ be the canonical bounded linear map.
For every sequence $(g_n)$ in $\mathcal F$, the
Eberlein--\v{S}mulian theorem gives a subsequence $(g_{n_k})$
converging weakly to some $g\in\mathcal F$. In particular,
$g_{n_k}(t)\to g(t)$ for every $t\in X$.
Since $\mathcal F$ is norm bounded and $\mu$ is finite,
the dominated convergence theorem gives
\[
\|Jg_{n_k}-Jg\|_2^2
=\int_X|g_{n_k}(t)-g(t)|^2\,d\mu(t)\longrightarrow0.
\]
Hence $J(\mathcal F)$ is norm compact.
The map $J$ is weak-to-weak continuous, and the weak and norm
topologies agree on the norm compact set $J(\mathcal F)$.
Therefore $J|_{\mathcal F}$ is weak-to-norm continuous.
\end{proof}

\begin{theorem}\label{commutative Bessel equality}
Let $(M,\Gamma)$ be a commutative Hopf-von Neumann algebra.
Then
\[
{\rm WAP}_{{\rm Bessel},\pi_u}(M_*)={\rm WAP}(M_*).
\]
\end{theorem}

\begin{proof}
By the Gelfand representation, there exists a compact Hausdorff
space $X$ such that $M=C(X)$. Let $x\in {\rm WAP}(M_*)$,
and let $(e_\lambda)_{\lambda\in\Lambda}$ be an orthonormal
basis for $H_u$.

For $t\in X$, denote evaluation at $t$ by $\delta_t\in M^*$,
and put
$a_t=R_x^{**}(\delta_t)\in C(X).$
For every $\nu\in M^*$, the pairing formula for $\Gamma_u(x)$
gives
\[
\nu(a_t)
=\nu\bigl(R_x^{**}(\delta_t)\bigr)
=\delta_t\bigl(L_x^{**}(\nu)\bigr).
\]
Since $L_x^{**}(\nu)\in C(X)$, the map $t\mapsto a_t$ is
weakly continuous. Consequently,
$\mathcal F=\{a_t:t\in X\}$ 
is a weakly compact subset of $C(X)$.

Fix $\zeta\in H_u$. Let $\mu_\zeta$ be the finite positive
regular Borel measure on $X$ representing the functional with the rule 
$ \mu_\zeta(a) = \langle\pi_u(a)\zeta,\zeta\rangle.$
Then
\[
\|\pi_u(a)\zeta\|^2
=\int_X|a(s)|^2\,d\mu_\zeta(s)
\qquad(a\in C(X)).
\]
By Lemma~\ref{DCT}, the map $t\mapsto a_t$ is norm continuous
as a map into $L^2(X,\mu_\zeta)$. Thus the function
\[
f(t)=\|\pi_u(a_t)\zeta\|^2
=\int_X|a_t(s)|^2\,d\mu_\zeta(s)
\]
belongs to $C(X)$.
For $\lambda\in\Lambda$, put
\[
b_\lambda
=(\omega_{\zeta,e_\lambda}\otimes {\rm id})(\Gamma_u(x))
=L_x^{**}(\omega_{\zeta,e_\lambda}),
\qquad f_\lambda=b_\lambda^*b_\lambda\in C(X)^+.
\]
For each $t\in X$, we have
\[
b_\lambda(t)
=\delta_t\bigl(L_x^{**}(\omega_{\zeta,e_\lambda})\bigr)
=\omega_{\zeta,e_\lambda}(a_t)
=\langle\pi_u(a_t)\zeta,e_\lambda\rangle.
\]
Hence, by Parseval's identity,
\begin{align*}
\sum_{\lambda\in F}f_\lambda(t)
=\sum_{\lambda\in F}
  |\langle\pi_u(a_t)\zeta,e_\lambda\rangle|^2
\longrightarrow\|\pi_u(a_t)\zeta\|^2=f(t)
\qquad(F\Subset\Lambda).
\end{align*}
Thus $\bigl(\sum_{\lambda\in F}f_\lambda\bigr)_{F\Subset\Lambda}$
is an increasing net of continuous positive functions whose
pointwise limit $f$ is continuous. Dini's theorem implies
\[
\left\|f-\sum_{\lambda\in F}f_\lambda\right\|_\infty
\longrightarrow0.
\]
In other words,
$\left[L_x^{**}(\omega_{\zeta,e_\lambda})
\right]_{\lambda\in\Lambda}\in K_{\Lambda,1}(M)$, for every $\zeta\in H_u$.
Since $x^*\in {\rm WAP}(M_*)$, applying the same argument to
$x^*$ and taking adjoints gives
\[
\left[L_x^{**}(\omega_{e_\lambda,\eta})
\right]_{\lambda\in\Lambda}\in K_{1,\Lambda}(M)
\qquad(\eta\in H_u).
\]
Interchanging $L_x^{**}$ and $R_x^{**}$ throughout the argument
gives the same column and row conditions for $R_x^{**}$.

Finally, let $(\eta_\alpha)_{\alpha\in I}$ be a Bessel family,
and let $S:\ell^2(I)\to H_u$ be its synthesis operator.
Identifying $H_u$ with $\ell^2(\Lambda)$, for
$T=L_x^{**}$ or $R_x^{**}$ we have
\[
\left[T(\omega_{\zeta,\eta_\alpha})\right]_{\alpha\in I}
=S^*\left[T(\omega_{\zeta,e_\lambda})\right]_{\lambda\in\Lambda}
\in K_{I,1}(M).
\]
Similarly, for a Bessel family $(\zeta_\alpha)_{\alpha\in I}$
with synthesis operator $S$, we have
\[
\left[T(\omega_{\zeta_\alpha,\eta})\right]_{\alpha\in I}
=\left[T(\omega_{e_\lambda,\eta})\right]_{\lambda\in\Lambda}S
\in K_{1,I}(M).
\]
Therefore $x\in {\rm WAP}_{{\rm Bessel},\pi_u}(M_*)$.
The reverse inclusion follows from the definition.
\end{proof}

Thus we arrive at the following corollary.

\begin{corollary}\label{commutative WAP algebra}
Let $(M,\Gamma)$ be a commutative Hopf-von Neumann algebra.
Then ${\rm WAP}(M_*)$ is a unital $C^*$-subalgebra of $M$.
\end{corollary}

\section{A coproduct characterization of weakly almost periodic functionals}

For locally compact quantum group $\mathbb{G}$, we denote
${\rm WAP}(M_*)$ by ${\rm WAP}(L^1(\mathbb{G}))$, where
$M=L^\infty(\mathbb{G})$. The theorem below characterizes this
space, for co-amenable $\mathbb{G}$, in terms of the
co-multiplication. 

\begin{theorem}\label{QM}
Let $\mathbb{G}$ be a co-amenable locally compact quantum group.
Then
\[
\begin{aligned}
{\rm WAP}(L^1(\mathbb{G}))
=\bigl\{x\in L^\infty(\mathbb{G}):\;&
\Gamma(x)\in QM(C_0(\mathbb{G})\otimes
                    {\rm WAP}(L^1(\mathbb{G})))\\
&\cap QM({\rm WAP}(L^1(\mathbb{G}))\otimes
                    C_0(\mathbb{G}))\bigr\}.
\end{aligned}
\]
\end{theorem}

\begin{proof}
Let $x\in {\rm WAP}(L^1(\mathbb{G}))$. Since this space is an
$L^1(\mathbb{G})$-bimodule, we have
$(\mu\otimes{\rm id})\Gamma(x)\in
{\rm WAP}(L^1(\mathbb{G}))$ for each $\mu\in L^1(\mathbb{G})$.
We show that
\[
(a\otimes1)\Gamma(x)(b\otimes1)
\in C_0(\mathbb{G})\otimes{\rm WAP}(L^1(\mathbb{G}))
\qquad(a,b\in C_0(\mathbb{G})).
\]
Fix $\mu\in L^1(\mathbb{G})$ and choose a normal extension
$\widetilde\mu\in B(L^2(\mathbb{G}))_*$ of $\mu$.
Since $C_0(\widehat{\mathbb{G}})$ acts non-degenerately on
$L^2(\mathbb{G})$, Cohen's factorization theorem gives
\[
\widetilde\mu=\hat b\cdot\mu_0\cdot\hat a,
\qquad
\hat a,\hat b\in C_0(\widehat{\mathbb{G}}),\quad
\mu_0\in B(L^2(\mathbb{G}))_*.
\]
Here the dual algebra is represented in its right regular
representation, as in the realization of $V$.
Using co-associativity and $\Gamma(x)=V(x\otimes1)V^*$, we obtain
\begin{align*}
&(a\otimes1)\Gamma((\mu\otimes{\rm id})\Gamma(x))(b\otimes1)\\
&=(\widetilde\mu\otimes{\rm id}\otimes{\rm id})
 \bigl((1\otimes a\otimes1)V_{12}V_{13}
 (x\otimes1\otimes1)V_{13}^*V_{12}^*
 (1\otimes b\otimes1)\bigr)\\
&=(\mu_0\otimes{\rm id}\otimes{\rm id})
 \bigl((\hat a\otimes a\otimes1)V_{12}V_{13}
 (x\otimes1\otimes1)V_{13}^*V_{12}^*
 (\hat b\otimes b\otimes1)\bigr).
\end{align*}
Since $V\in M(C_0(\widehat{\mathbb{G}})\otimes C_0(\mathbb{G}))$,
the elements $(\hat a\otimes a)V$ and $V^*(\hat b\otimes b)$
belong to $C_0(\widehat{\mathbb{G}})\otimes C_0(\mathbb{G})$.
Approximating them in norm by finite sums
$\sum_j\hat a_j\otimes a_j$ and
$\sum_k\hat b_k\otimes b_k$, respectively, shows that the last
expression is a norm limit of sums of the form
\[
\sum_{j,k}a_jb_k\otimes
 \bigl((\mu_{jk}\otimes{\rm id})\Gamma(x)\bigr),
\qquad
\mu_{jk}=
 (\hat b_k\cdot\mu_0\cdot\hat a_j)|_{L^\infty(\mathbb{G})}
 \in L^1(\mathbb{G}).
\]
Consequently,
\[
(a\otimes1)\Gamma((\mu\otimes{\rm id})\Gamma(x))(b\otimes1)
\in C_0(\mathbb{G})\otimes{\rm WAP}(L^1(\mathbb{G})).
\]

Let $(\epsilon_\alpha)$ be a bounded approximate identity of
$L^1(\mathbb{G})$ consisting of normal states; such an
approximate identity exists by co-amenability
\cite[Section 3]{Runde3}. Then
$(\epsilon_\alpha\otimes{\rm id})\Gamma(x)\to x$ in the
weak$^*$ topology. Since $L_x$ is weakly compact, this net is
relatively weakly compact and all its weak cluster points equal
$x$. Hence it converges weakly to $x$, and therefore
\[
\begin{aligned}
&(a\otimes1)
 \Gamma((\epsilon_\alpha\otimes{\rm id})\Gamma(x))(b\otimes1)
 \xrightarrow{weak}(a\otimes1)\Gamma(x)(b\otimes1)\\
&\hspace{35mm}\in
 C_0(\mathbb{G})\otimes{\rm WAP}(L^1(\mathbb{G})),
\end{aligned}
\]
because this tensor product is a norm-closed linear subspace
and hence is weakly closed. The corresponding argument with
the left fundamental unitary $W$, using
$\Gamma(x)=W^*(1\otimes x)W$, gives
\[
(1\otimes a)\Gamma(x)(1\otimes b)
\in {\rm WAP}(L^1(\mathbb{G}))\otimes C_0(\mathbb{G}).
\]
This proves the required inclusion.

Conversely, suppose that $x\in L^\infty(\mathbb{G})$ and
$\Gamma(x)\in QM(C_0(\mathbb{G})\otimes
{\rm WAP}(L^1(\mathbb{G})))$.
Passing to a subnet, we may assume that
$\epsilon_\alpha|_{C_0(\mathbb{G})}$ converges in
$\sigma(\mathcal M(\mathbb{G}),C_0(\mathbb{G}))$ to the bounded
counit $\epsilon$; see \cite[Proposition 5.2 and its proof]{Runde3}.
Choose a positive contraction $a\in C_0(\mathbb{G})$ such that
$\epsilon(a)=1$. Since $\epsilon$ is multiplicative,
\[
\epsilon_\alpha((1-a)^2)
=1-2\epsilon_\alpha(a)+\epsilon_\alpha(a^2)\longrightarrow0.
\]
The Cauchy-Schwarz inequality for states gives
\[
\|a\cdot\epsilon_\alpha\cdot a-\epsilon_\alpha\|
\leq 2\epsilon_\alpha((1-a)^2)^{1/2}\longrightarrow0.
\]
On the other hand,
\[
Z=(a\otimes1)\Gamma(x)(a\otimes1)
\in C_0(\mathbb{G})\otimes{\rm WAP}(L^1(\mathbb{G})).
\]
Approximating $Z$ in norm by finite sums of elementary tensors,
we obtain
\[
\begin{aligned}
(a\cdot\epsilon_\alpha\cdot a\otimes{\rm id})\Gamma(x)
&=(\epsilon_\alpha\otimes{\rm id})(Z)\\
&\longrightarrow(\epsilon\otimes{\rm id})(Z)
\in{\rm WAP}(L^1(\mathbb{G}))
\end{aligned}
\]
in norm. The same net converges weak$^*$ to $x$, because
$(\epsilon_\alpha\otimes{\rm id})\Gamma(x)\to x$ weak$^*$
and $\|a\cdot\epsilon_\alpha\cdot a-\epsilon_\alpha\|\to0$.
Thus $x=(\epsilon\otimes{\rm id})(Z)$ belongs to
${\rm WAP}(L^1(\mathbb{G}))$, as required.
\end{proof}

\begin{remark}\label{regular WAP multipliers}
A multiplicative unitary $V$ is regular if
\[
\overline{\{({\rm id}\otimes\omega)(\Sigma V):
             \omega\in B(L^2(\mathbb{G}))_*\}}^{\|\cdot\|}
=K(L^2(\mathbb{G})).
\]
For a regular multiplicative unitary, we have
\[
C_0(\widehat{\mathbb{G}})\otimes C_0(\mathbb{G})
=\overline{\operatorname{span}}^{\|\cdot\|}
 \{(\hat a\otimes1)V(1\otimes a):
   \hat a\in C_0(\widehat{\mathbb{G}}),\ a\in C_0(\mathbb{G})\};
\]
see \cite[Proposition 3.6(d)]{Baaj}.
Following the use of regularity in \cite{Salm}, we can replace
quasi-multipliers by multipliers in Theorem~\ref{QM}.
Indeed, if $\mathbb{G}$ is co-amenable and regular, then
\begin{equation}\label{Mcase}
\begin{aligned}
{\rm WAP}(L^1(\mathbb{G}))
=\bigl\{x\in L^\infty(\mathbb{G}):\;&
\Gamma(x)\in M(C_0(\mathbb{G})\otimes
                   {\rm WAP}(L^1(\mathbb{G})))\\
&\cap M({\rm WAP}(L^1(\mathbb{G}))\otimes
                   C_0(\mathbb{G}))\bigr\}.
\end{aligned}
\end{equation}

To see this, let $x\in{\rm WAP}(L^1(\mathbb{G}))$,
$a\in C_0(\mathbb{G})$ and $\mu\in L^1(\mathbb{G})$.
As in the preceding proof, choose a normal extension
$\widetilde\mu=\hat a_2\cdot\mu_0\cdot\hat a_1$, where
$\hat a_1,\hat a_2\in C_0(\widehat{\mathbb{G}})$ and
$\mu_0\in B(L^2(\mathbb{G}))_*$. Then
\begin{align*}
\Gamma((\mu\otimes{\rm id})\Gamma(x))(a\otimes1)
=(\mu_0\otimes{\rm id}\otimes{\rm id})
 \bigl((\hat a_1\otimes1\otimes1)V_{12}V_{13}
 (x\otimes1\otimes1)V_{13}^*V_{12}^*
 (\hat a_2\otimes a\otimes1)\bigr).
\end{align*}
Approximate $V^*(\hat a_2\otimes a)$ in norm by finite sums
$\sum_j\hat b_j\otimes b_j$. For each $j$, regularity gives
\[
(\hat a_1\otimes1)V(1\otimes b_j)
\in C_0(\widehat{\mathbb{G}})\otimes C_0(\mathbb{G}),
\]
so this element can in turn be approximated in norm by finite
sums $\sum_k\hat a_{jk}\otimes b_{jk}$.
Since $1\otimes b_j\otimes1$ commutes with
$V_{13}(x\otimes1\otimes1)V_{13}^*$, the displayed slice is
a norm limit of sums of the form
\[
\sum_{j,k}b_{jk}\otimes
 \bigl((\mu_{jk}\otimes{\rm id})\Gamma(x)\bigr),
\qquad
\mu_{jk}=
 (\hat b_j\cdot\mu_0\cdot\hat a_{jk})|_{L^\infty(\mathbb{G})}.
\]
Thus
\[
\Gamma((\mu\otimes{\rm id})\Gamma(x))(a\otimes1)
\in C_0(\mathbb{G})\otimes{\rm WAP}(L^1(\mathbb{G})).
\]
Taking $\mu=\epsilon_\alpha$ and passing to the weak limit as
in Theorem~\ref{QM} gives the same inclusion for
$\Gamma(x)(a\otimes1)$. Applying this to $x^*$ and taking
adjoints gives the inclusion for $(a\otimes1)\Gamma(x)$.
The corresponding argument with $W$ gives the two inclusions
on the second tensor factor. The reverse inclusion in
\eqref{Mcase} follows from Theorem~\ref{QM}.

Another condition which gives \eqref{Mcase} is the existence
of a bounded approximate identity $(e_\beta)$ for
$C_0(\mathbb{G})$ in its center; compare \cite{Runde3}.
Indeed, since $C_0(\mathbb{G})$ is weak$^*$ dense in
$L^\infty(\mathbb{G})$, each $e_\beta$ is also central in
$L^\infty(\mathbb{G})$. Hence, for
$x\in{\rm WAP}(L^1(\mathbb{G}))$ and $a\in C_0(\mathbb{G})$,
\[
\begin{aligned}
\Gamma(x)(a\otimes1)
=\lim_\beta(e_\beta\otimes1)\Gamma(x)(a\otimes1)
\in C_0(\mathbb{G})\otimes{\rm WAP}(L^1(\mathbb{G})),
\end{aligned}
\]
where the limit is in norm. The other three multiplier
conditions follow in the same way from Theorem~\ref{QM}.

Finally, the same arguments apply to the Bessel algebra
${\rm WAP}_{{\rm Bessel},\pi_u}(L^1(\mathbb{G}))$.
For completeness, this algebra is an $L^1(\mathbb{G})$-bimodule.
Indeed, for $\mu\in L^1(\mathbb{G})$ put
$S_\mu(a)=(\mu\otimes{\rm id})\Gamma(a)$.
Co-associativity gives, for $y=S_\mu(x)$,
\[
L_y^{**}=L_x^{**}\circ S_\mu^*,
\qquad
R_y^{**}=S_\mu\circ R_x^{**},
\]
where the ranges of the biadjoints are identified with their
canonical copies in $L^\infty(\mathbb{G})$.
The maps $S_\mu$ and $S_\mu^*$ are completely bounded.
Thus $S_\mu^*$ preserves the bounded columns and rows over
$L^\infty(\mathbb{G})^*$, identified earlier with
$Q_c^{H_u}$ and $Q_r^{H_u}$, while $S_\mu$ preserves the
corresponding $K$-columns and $K$-rows.
This proves the Bessel conditions for $y$; the other module
action is treated in the same way. Therefore, under either
of the preceding additional hypotheses,
\[
\begin{aligned}
{\rm WAP}_{{\rm Bessel},\pi_u}(L^1(\mathbb{G}))
=\bigl\{x\in L^\infty(\mathbb{G}):
 \Gamma(x)\in
 & M(C_0(\mathbb{G})\otimes
       {\rm WAP}_{{\rm Bessel},\pi_u}(L^1(\mathbb{G})))\\
\hspace{34mm}\cap
 & M({\rm WAP}_{{\rm Bessel},\pi_u}(L^1(\mathbb{G}))
       \otimes C_0(\mathbb{G}))\bigr\}.
\end{aligned}
\]
Here co-amenability is retained in both cases.
\end{remark}

\section{Left uniform continuity on locally compact quantum groups}

In this section, we study left uniformly continuous elements of
locally compact quantum groups. Since the left uniformly continuous
functions on a locally compact group form a $C^*$-algebra, it is
natural to ask under what conditions on $\mathbb{G}$ the space
${\rm LUC}(L^1(\mathbb{G}))$ is also a $C^*$-algebra. This is known
when $\mathbb{G}$ is semi-regular, and also when $\mathbb{G}$ is
co-amenable and $C_0(\mathbb{G})$ has a bounded approximate
identity in its center \cite{HuN,Runde3}.
We first establish the $C^*$-algebra property under a weaker condition involving a single central element.
We then characterize ${\rm LUC}(\mathbb{G})$ as the greatest
left introverted $C^*$-subalgebra of
${\rm LUC}(L^1(\mathbb{G}))$ satisfying the left Bessel condition
introduced below.

\begin{proposition}\label{LUC QM}
Let $\mathbb{G}$ be a co-amenable locally compact quantum group.
Then
\[
\begin{aligned}
{\rm LUC}(L^1(\mathbb{G}))
=\bigl\{x\in C_b(\mathbb{G}):\;&\Gamma(x)\in
QM\bigl(C_0(\mathbb{G})\otimes
              {\rm LUC}(L^1(\mathbb{G}))\bigr)\bigr\}.
\end{aligned}
\]
\end{proposition}

\begin{proof}
For $x\in{\rm LUC}(L^1(\mathbb{G}))$ and
$\mu\in L^1(\mathbb{G})$, we have
$(\mu\otimes{\rm id})\Gamma(x)\in
{\rm LUC}(L^1(\mathbb{G}))$. The finite-sum approximation in
the proof of Theorem~\ref{QM}, which uses only this invariance,
therefore gives
\[
(a\otimes1)\Gamma(x\star\mu)(b\otimes1)
\in C_0(\mathbb{G})\otimes{\rm LUC}(L^1(\mathbb{G}))
\qquad(a,b\in C_0(\mathbb{G})).
\]
If $(\epsilon_\alpha)$ is a bounded approximate identity of
$L^1(\mathbb{G})$, then $x\star\epsilon_\alpha\to x$ in norm.
Taking $\mu=\epsilon_\alpha$ proves the required inclusion.

Conversely, let $x\in C_b(\mathbb{G})$ satisfy the stated
quasi-multiplier condition, and let $\epsilon$ be the bounded
counit of $C_0(\mathbb{G})$. Choose $a\in C_0(\mathbb{G})$
such that $\epsilon(a)=1$. The counit identity extends strictly
to $C_b(\mathbb{G})$, and hence
\[
x=(\epsilon\otimes{\rm id})
       \bigl((a\otimes1)\Gamma(x)(a\otimes1)\bigr)
\in{\rm LUC}(L^1(\mathbb{G})).
\]
\end{proof}

The following result uses a single central element on which the
counit is non-zero, in place of the central approximate identity
assumption in \cite{Runde3}.

\begin{theorem}\label{central LUC algebra}
Let $\mathbb{G}$ be a co-amenable locally compact quantum group,
and let $\epsilon\in\mathcal M(\mathbb{G})$ be its bounded
counit. If there is an element $a$ in the center of
$C_0(\mathbb{G})$ such that $\epsilon(a)\neq0$, then
${\rm LUC}(L^1(\mathbb{G}))$ is a $C^*$-algebra.
\end{theorem}

\begin{proof}
Replacing $a$ by a suitable continuous function of $a^*a$, we
may assume that $0\leq a\leq1$ and $\epsilon(a)=1$.
Since $C_0(\mathbb{G})$ is weak$^*$ dense in
$L^\infty(\mathbb{G})$, the element $a$ is also central in
$L^\infty(\mathbb{G})$.
Choose a bounded approximate identity $(\epsilon_\alpha)$ of
$L^1(\mathbb{G})$ consisting of normal states, and pass to a
subnet so that
\[
\epsilon_\alpha|_{C_0(\mathbb{G})}
\longrightarrow\epsilon
\quad\hbox{in }\sigma(\mathcal M(\mathbb{G}),C_0(\mathbb{G})).
\]

Let $x,z\in{\rm LUC}(L^1(\mathbb{G}))$. By Cohen's
factorization theorem, we may write $x=y\star\mu$, where
$y\in{\rm LUC}(L^1(\mathbb{G}))$ and
$\mu\in L^1(\mathbb{G})$. Factor once more, using the
$C_0(\mathbb{G})$-bimodule structure of $L^1(\mathbb{G})$,
to obtain
$\mu=e'\cdot\widetilde\mu\cdot e,$ for some $ e,e'\in C_0(\mathbb{G}),$  and
$\widetilde\mu\in L^1(\mathbb{G}).$
Put
\[
P=(e\otimes1)\Gamma(y)(e'\otimes1),
\qquad Q=(a\otimes1)\Gamma(z)(a\otimes1).
\]
By Proposition~\ref{LUC QM}, both $P$ and $Q$ belong to
$C_0(\mathbb{G})\otimes{\rm LUC}(L^1(\mathbb{G}))$.
Using centrality of $a$ and co-associativity, we obtain
\begin{align*}
&(xz)\star(a\cdot\epsilon_\alpha\cdot a^2)\\
&=(\mu\otimes a\cdot\epsilon_\alpha\cdot a^2
                         \otimes{\rm id})
 \bigl(( {\rm id}\otimes\Gamma)\Gamma(y)
                          (1\otimes\Gamma(z))\bigr)\\
&=(\mu\otimes\epsilon_\alpha\cdot a\otimes{\rm id})
 \bigl(( {\rm id}\otimes\Gamma)\Gamma(y)(1\otimes Q)\bigr)\\
&=(\widetilde\mu\otimes\epsilon_\alpha\otimes{\rm id})
 \bigl((1\otimes a\otimes1)
       ({\rm id}\otimes\Gamma)(P)(1\otimes Q)\bigr).
\end{align*}
Approximate $P$ and $Q$ in norm by finite sums
$\sum_j f_j\otimes g_j$ and $\sum_k c_k\otimes d_k$,
respectively, where $f_j,c_k\in C_0(\mathbb{G})$ and
$g_j,d_k\in{\rm LUC}(L^1(\mathbb{G}))$.
The last expression is then approximated, uniformly in $\alpha$,
by sums of the form
\[
\sum_{j,k}\widetilde\mu(f_j)
 \bigl[(\epsilon_\alpha\otimes{\rm id})
       ((a\otimes1)\Gamma(g_j)(c_k\otimes1))\bigr]d_k.
\]
For each $j,k$, Proposition~\ref{LUC QM} gives
\[
(a\otimes1)\Gamma(g_j)(c_k\otimes1)
\in C_0(\mathbb{G})\otimes{\rm LUC}(L^1(\mathbb{G})).
\]
Consequently, its slices by $\epsilon_\alpha$ converge in norm
to its slice by $\epsilon$. The uniformity of the preceding
approximations follows from $\|\epsilon_\alpha\|=1$ and the
boundedness of the slice maps. It follows that
$((xz)\star(a\cdot\epsilon_\alpha\cdot a^2))_\alpha$ is a
norm-Cauchy net.

Since $a$ is central and $0\leq a\leq1$, we also have
\[
\|a\cdot\epsilon_\alpha\cdot a^2-\epsilon_\alpha\|
=1-\epsilon_\alpha(a^3)\longrightarrow0.
\]
Thus $(xz)\star(a\cdot\epsilon_\alpha\cdot a^2)$ converges
weak$^*$ to $xz$, and its norm limit must be $xz$.
Each term belongs to ${\rm LUC}(L^1(\mathbb{G}))$, which is
norm closed. Hence $xz\in{\rm LUC}(L^1(\mathbb{G}))$.
As this space is already a self-adjoint operator system, it is
a $C^*$-algebra.
\end{proof}

We next record the quasi-multiplier version of
\cite[Lemma A.4]{KV1}. Functionals on $A$ are evaluated on
$A^{**}$ through their canonical normal extensions.

\begin{theorem}\label{stconv}
Let $A$ be a $C^*$-algebra, let $(a_i)_{i\in I}$ be an
increasing net in $QM(A)^+$, and let $a\in QM(A)^+$ satisfy
\[
\omega(a)=\sup_{i\in I}\omega(a_i)
\qquad(\omega\in A^*_+).
\]
Then $(a_i)$ converges strictly to $a$, in the sense that
\[
\|(a-a_i)b\|+\|b(a-a_i)\|\longrightarrow0
\qquad(b\in A).
\]
\end{theorem}

\begin{proof}
For $b\in A$, the elements $b^*a_i b$ and $b^*ab$ belong to
$A^+$. Applying the hypothesis to
$t\mapsto\omega(b^*tb)$, for $\omega\in A^*_+$, and then
\cite[Lemma A.3]{KV1}, we obtain
\[
\|b^*(a-a_i)b\|\longrightarrow0.
\]
Since $0\leq a-a_i\leq a$, it follows that
\[
\|(a-a_i)b\|^2
=\|b^*(a-a_i)^2b\|
\leq\|a\|\,\|b^*(a-a_i)b\|\longrightarrow0.
\]
Replacing $b$ by $b^*$ and taking adjoints proves the other
convergence.
\end{proof}

We also need the following slice expansion.

\begin{proposition}[{\cite[Lemma A.5]{KV1}}]\label{stst}
Let $A$ and $B$ be $C^*$-algebras, let $\pi:A\to B(H)$ be
a non-degenerate $*$-representation, and let $(e_i)_{i\in I}$
be an orthonormal basis for $H$. If $x,y\in M(A\otimes B)$
and $\zeta,\eta\in H$, then the net
\[
\left(
\sum_{i\in F}
 (\omega_{\eta,e_i}\otimes{\rm id})(y)^*
 (\omega_{\zeta,e_i}\otimes{\rm id})(x)
\right)_{F\Subset I}
\]
is bounded and converges strictly in $M(B)$ to
$(\omega_{\zeta,\eta}\otimes{\rm id})(y^*x)$.
Here the vector functionals on $A$ are taken with respect
to $\pi$. If $B$ is unital, the convergence is in norm.
\end{proposition}

Fix an orthonormal basis $(e_i)_{i\in I}$ for
$L^2(\mathbb{G})$. We say that a left invariant
$C^*$-subalgebra $C$ of $L^\infty(\mathbb{G})$ satisfies the
\emph{left Bessel condition} if
\begin{equation}\tag{$\mathrm B_L$}\label{left Bessel condition}
\left[(\omega_{\zeta,e_i}\otimes{\rm id})\Gamma(x)\right]_{i\in I}
\in K_{I,1}(C)
\qquad(x\in C,\ \zeta\in L^2(\mathbb{G})).
\end{equation}
Equivalently, the net
\[
\sum_{i\in F}
 \bigl((\omega_{\zeta,e_i}\otimes{\rm id})\Gamma(x)\bigr)^*
 \bigl((\omega_{\zeta,e_i}\otimes{\rm id})\Gamma(x)\bigr).
\]
converges in norm as $F$ runs through the finite subsets of $I$.
This condition is independent of the chosen orthonormal basis.
Indeed, if $(\eta_j)_{j\in J}$ is a Bessel family and
$T:\ell^2(J)\to L^2(\mathbb{G})$ is its synthesis operator,
$T\delta_j=\eta_j$, then
\[
\left[(\omega_{\zeta,\eta_j}\otimes{\rm id})\Gamma(x)\right]_{j\in J}
=T^*
 \left[(\omega_{\zeta,e_i}\otimes{\rm id})\Gamma(x)\right]_{i\in I},
\]
where $L^2(\mathbb{G})$ is identified with $\ell^2(I)$ through
the chosen basis. Multiplication by $T^*$ preserves compact
columns. Thus \eqref{left Bessel condition} is equivalent to
the same condition for all Bessel families on $L^2(\mathbb{G})$.
Only normal vector functionals are involved here; this is
distinct from the universal-representation condition defining
${\rm WAP}_{{\rm Bessel},\pi_u}$.

\begin{proposition}\label{bigluc}
Let $\mathbb{G}$ be a co-amenable locally compact quantum group.
Then
\[
{\rm LUC}(\mathbb{G})
=\bigl\{x\in C_b(\mathbb{G}):
 \Gamma(x)\in M(C_0(\mathbb{G})\otimes C_b(\mathbb{G}))\bigr\}
\]
is the greatest left introverted $C^*$-subalgebra of
${\rm LUC}(L^1(\mathbb{G}))$ satisfying
\eqref{left Bessel condition}.
\end{proposition}

\begin{proof}
Let $C$ be a left introverted $C^*$-subalgebra of
${\rm LUC}(L^1(\mathbb{G}))$ satisfying
\eqref{left Bessel condition}, and let $y\in C$.
Since $C$ is left invariant and
$y\star\epsilon_\alpha\to y$ in norm, Cohen's factorization
theorem gives
\[
y=x\star\mu,
\qquad x\in C,\quad\mu\in L^1(\mathbb{G}).
\]
Let $a\in C_0(\mathbb{G})$. We show that
$\Gamma(y)(a\otimes1)\in C_0(\mathbb{G})\otimes C_b(\mathbb{G})$.

As $L^\infty(\mathbb{G})$ acts in standard form on
$L^2(\mathbb{G})$, write
$\mu=\omega_{\xi,\eta}|_{L^\infty(\mathbb{G})}$, with
$\xi,\eta\in L^2(\mathbb{G})$.
The non-degenerate action of $C_0(\widehat{\mathbb{G}})$ on
$L^2(\mathbb{G})$ allows us to factor
$\xi=\hat a\xi'$, where
$\hat a\in C_0(\widehat{\mathbb{G}})$ and
$\xi'\in L^2(\mathbb{G})$.
Put $\mu'=\omega_{\xi',\eta}\in B(L^2(\mathbb{G}))_*$.
Thus $\hat a\cdot\mu'$ restricts to $\mu$ on
$L^\infty(\mathbb{G})$. As before, the dual algebra is in its
right regular representation. We obtain
\begin{align*}
\Gamma(y)(a\otimes1)
&=(\mu\otimes{\rm id}\otimes{\rm id})
       (({\rm id}\otimes\Gamma)\Gamma(x))(a\otimes1)\\
&=(\mu'\otimes{\rm id}\otimes{\rm id})
 \bigl(V_{12}V_{13}(x\otimes1\otimes1)V_{13}^*V_{12}^*
                   (\hat a\otimes a\otimes1)\bigr).
\end{align*}
Since
\[
V^*(\hat a\otimes a)
\in C_0(\widehat{\mathbb{G}})\otimes C_0(\mathbb{G}),
\]
we may approximate it in norm by finite sums
$\sum_j\hat b_j\otimes b_j$.
It is therefore enough to show that, for
$\zeta,\eta\in L^2(\mathbb{G})$ and $b\in C_0(\mathbb{G})$,
\[
(\omega_{\zeta,\eta}\otimes{\rm id}\otimes{\rm id})
 \bigl(V_{12}V_{13}(x\otimes1\otimes1)V_{13}^*\bigr)
 (b\otimes1)
\in C_0(\mathbb{G})\otimes C_b(\mathbb{G}).
\]
Indeed, the terms arising from the finite sums have this form
with $\zeta=\hat b_j\xi'$ and $b=b_j$.

Put
\[
u_i=(\omega_{e_i,\eta}\otimes{\rm id})(V),
\qquad
d_i=(\omega_{\zeta,e_i}\otimes{\rm id})\Gamma(x).
\]
Then $u_i\in C_0(\mathbb{G})$, $d_i\in C$, and the usual
matrix multiplication identity gives
\begin{align*}
(\omega_{\zeta,\eta}\otimes{\rm id}\otimes{\rm id})
 \bigl(V_{12}V_{13}(x\otimes1\otimes1)V_{13}^*\bigr)
 (b\otimes1)
=\sum_{i\in I}u_i b\otimes d_i,
\end{align*}
initially with ultraweak convergence. In the row and column
notation used above, put
\[
X=[u_i\otimes1]_{i\in I},
\qquad
Y=[b\otimes d_i]_{i\in I}^{\rm t}.
\]
For $F\Subset I$, let $X_F$ and $Y_F$ be obtained by replacing
the entries outside $F$ by zero. Thus
\[
X_FY_F=\sum_{i\in F}u_i b\otimes d_i
\in C_0(\mathbb{G})\otimes C.
\]
Unitarity of $V$ gives
\[
\sum_{i\in F}u_i u_i^*\leq\|\eta\|^2 1,
\qquad \|X\|\leq\|\eta\|.
\]
In fact, Proposition~\ref{stst}, applied to $V^*$, shows that
these sums converge strictly to $\|\eta\|^2 1$.

For the column, put
\[
q=(\omega_{\zeta,\zeta}\otimes{\rm id})\Gamma(x^*x),
\qquad
z_F=b^*b\otimes\sum_{i\in F}d_i^*d_i,
\qquad z=b^*b\otimes q.
\]
The sums $\sum_{i\in F}d_i^*d_i$ converge ultraweakly to $q$.
By \eqref{left Bessel condition}, they also converge in norm,
so their norm limit is $q$. Since $x^*x\in C$ and $C$ is
left invariant, $q\in C$. Consequently,
\[
\|z-z_F\|\longrightarrow0,
\qquad
\|Y-Y_F\|^2=\|z-z_F\|\longrightarrow0.
\]
It follows that
\[
\|XY-X_FY_F\|
\leq\|X\|\,\|Y-Y_F\|
\leq\|\eta\|\,\|z-z_F\|^{1/2}\longrightarrow0.
\]
Thus the series $\sum_i u_i b\otimes d_i$ converges in norm
in $C_0(\mathbb{G})\otimes C$, and its norm limit equals
the displayed ultraweak limit. Returning to the finite-sum
approximations of $V^*(\hat a\otimes a)$, we conclude that
\[
\Gamma(y)(a\otimes1)
\in C_0(\mathbb{G})\otimes C_b(\mathbb{G}).
\]
Applying the same argument to $y^*$ and taking adjoints gives
$(a\otimes1)\Gamma(y)\in C_0(\mathbb{G})\otimes C_b(\mathbb{G})$.
Since $C_b(\mathbb{G})$ is unital, these two inclusions imply
$\Gamma(y)\in M(C_0(\mathbb{G})\otimes C_b(\mathbb{G})).$
Hence $C\subseteq{\rm LUC}(\mathbb{G})$.

It remains to check that ${\rm LUC}(\mathbb{G})$ itself belongs
to the stated class. Its unital $C^*$-algebra structure and
left introversion are established in \cite{Salm}. In the
co-amenable case, its inclusion in
${\rm LUC}(L^1(\mathbb{G}))$ follows from the quasi-multiplier
characterization in \cite{Runde3}.
For $x\in{\rm LUC}(\mathbb{G})$ and
$\zeta\in L^2(\mathbb{G})$, apply Proposition~\ref{stst} with
$A=C_0(\mathbb{G}),$ $B=C_b(\mathbb{G}),$
and with both multiplier arguments equal to $\Gamma(x)$.
Writing $d_i=(\omega_{\zeta,e_i}\otimes{\rm id})\Gamma(x)$,
we obtain
\[
\sum_{i\in F}d_i^*d_i
\xrightarrow{\rm strict}
(\omega_{\zeta,\zeta}\otimes{\rm id})\Gamma(x^*x).
\]
Because $C_b(\mathbb{G})$ is unital, this convergence is in
norm. Moreover, left invariance gives
$d_i\in{\rm LUC}(\mathbb{G})$. Since this algebra is norm
closed,
$[d_i]_{i\in I}\in K_{I,1}({\rm LUC}(\mathbb{G})).$
Thus ${\rm LUC}(\mathbb{G})$ satisfies
\eqref{left Bessel condition}. Together with the inclusion
proved above, this establishes the result.
\end{proof}
\begin{corollary}
    Let $\mathbb{G}$ be a co-amenable locally compact quantum group,
and let $\epsilon\in\mathcal M(\mathbb{G})$ be its bounded
counit. If there is an element $a$ in the center of
$C_0(\mathbb{G})$ such that $\epsilon(a)\neq0$, then
${\rm LUC}(L^1(\mathbb{G})) = {\rm LUC}(\mathbb{G})$
\end{corollary}

\section*{acknowledgment} 

This work was mainly carried out while the author was
affiliated with Ferdowsi University of Mashhad, Mashhad, Iran,
and was completed at the University of Oulu, Oulu, Finland. 

During the preparation of this manuscript, ChatGPT 6 was used to enhance the language and readability of the text. Following the use of this tool, the authors carefully reviewed and revised the content as necessary and take full responsibility for this version of the manuscript.

\end{document}